\input amstex
\font\tc=cmr10 \font\tcb=cmbx10 \font\tci=cmti10 \font\tcr=cmss10
\font\ttc=cmr8 \font\ttci=cmti8

\def\codim{\text{\rm codim}\,}

\define\pd#1#2{\dfrac{\partial#1}{\partial#2}}
\define\p2d#1#2#3{\dfrac{\partial^2#1}{\partial#2\partial#3}}
\NoBlackBoxes

\documentstyle{amsppt}
\baselineskip=2\baselineskip
\topmatter
\title{$L$-convex-concave sets in real projective space and $L$-duality
\footnotemark"*"}
\endtitle
\rightheadtext{$L$-convex-concave sets in real projective space and
 $L$-duality}
\footnotetext"*"{\ttc Khovanskii's work is partially supported by
Canadian Grant {\rm N~0GP0156833}. Novikov's work was supported by
the Killam grant of P. Milman and by James S. McDonnell
Foundation.}

\author{A.~Khovanskii,  D.~Novikov}
\endauthor

\address1. {Department of Mathematics, Toronto University, Toronto,
Canada}
\endaddress

\email{askold\@math.toronto.edu}
\endemail

\address2.
{Department of Mathematics, Purdue University, West Lafayette IN}
\endaddress
\email{dmitry\@math.purdue.edu}
\endemail

\abstract\nofrills{\ttc We define a class of $L$-convex-concave
subsets of $\Bbb{R}P^n$, where $L$ is a projective subspace of
dimension $l$ in $\Bbb{R}P^n$. These are sets whose sections by
any $(l+1)$-dimensional space $L'$ containing $L$ are convex and
concavely depend on $L'$. We introduce an $L$-duality for these
sets, and prove that the $L$-dual to an $L$-convex-concave set is
an $L^*$-convex-concave subset of $(\Bbb RP^n)^*$.  We discuss a
version of Arnold hypothesis for these sets and prove that it is
true (or wrong) for an $L$-convex-concave set and its $L$-dual
simultaneously. }
\endabstract
\endtopmatter
\document

\head{\tcb Introduction}
\endhead

\subhead{\tcb Convex-concave sets and Arnold hypothesis}
\endsubhead
{\tc The notion of  convexity is  usually defined for
subsets of affine spaces, but it can be  generalized for subsets
of projective spaces. Namely, a subset of  a projective space $\Bbb
RP^n$ is  called {\tci convex} if it doesn't intersect some hyperplane
$L\subset \Bbb RP^n$ and is convex in the affine space $\Bbb
RP^n\setminus L$. In  the very definition of the convex subset of a
projective space appears a hyperplane $L$. In  projective space
there are subspaces $L$ of different dimensions, not only
hyperplanes. For any subspace $L$ one can define a class of
$L$-convex-concave sets. These sets are the main object of
investigation in this paper. If $L$ is a
hyperplane then  this class coincides with the class of closed
convex sets lying in the  affine chart $\Bbb RP^n\setminus L$.
Here is the definition of $L$-convex-concave  sets.

A closed set $A\subset\Bbb RP^n$ is $L$-{\tci convex-concave} if:
1) the set $A$ doesn't intersect the projective subspace $L$, 2)
for any $(\dim L+1)$-dimensional subspace $N\subset\Bbb RP^n$
containing $L$ the section $A\cap N$ of the set $A$ by $N$ is
convex, 3) for any $(\dim L-1)$-dimensional subspace $T\subset L$
the complement to the projection of the set $A$ from the center
$T$ on the factor-space $\Bbb RP^n/T$ is an open convex set. }

\example{\tcr Example}{\tc In a projective space $\Bbb RP^n$ with
homogeneous coordinates $x_0:\dots :x_n$ one can consider a set
$A\subset\Bbb RP^n$ defined by the inequality $\{K(x)\leq 0\}$,
where $K$ is a non-degenerate quadratic form on $\Bbb R^{n+1}$.
Suppose that $K$ is positively defined on some $(k+1)$-dimensional
subspace, and is negatively defined on some $(n-k)$-dimensional
subspace. In other words, suppose that (up to a linear change of
coordinates) the form $K$ is of the form
$K(x)=x_0^2+\dots+x_k^2-x_{k+1}^2-\dots-x_n^2$. In this case the
set $A$ is $L$-convex-concave with respect to projectivization $L$ of any
$(k+1)$-dimensional subspace of $\Bbb R^{n+1}$ on which $K$ is
positively defined.}
\endexample

{\tc We are mainly interested in the  following hypothesis.}

\proclaim{\tcr The Main Hypothesis}{ \tci Any $L$-convex-concave
subset $A$ of an  $n$-dimensional projective space contains a
projective subspace $M$ of dimension equal to $(n-1-\dim L)$.}
\endproclaim

{\tc Note that any projective subspace of dimension bigger than
$(n-1-\dim L)$ necessarily intersects $L$, so it cannot be
contained in $A$. For the quadratic set $A$ from the previous
example the Main Hypothesis is evidently true: as $M$ one can take
projectivization of any $(n-k)$-dimensional subspaces of $\Bbb
R^{n+1}$ on which $K$  is negatively defined.

For an $L$-convex-concave  set $A$ with a smooth non-degenerate
boundary $B$ the Main Hypothesis is a particular case of the following
hypothesis due to Arnold, see \cite{Ar1, Ar2}.}

\proclaim{\tcr Arnold hypothesis }{ \tci  Let $B\subset\Bbb
RP^n$ be  a connected
 smooth hypersurface  bounding some domain $U\subset \Bbb RP^n$.
Suppose that at any point of $B$ the  second fundamental form of
$B$ with respect to the outward normal vector is nondegenerate.
Suppose that this form has a (necessarily constant) signature
$(n-k-1,k)$, i.e. at each point $b\in B$ the  restriction of the
second quadratic form to some $k$-dimensional subspace of $T_b B$
is negatively defined and its restriction to some
$(n-k-1)$-dimensional subspace of $T_b B$ is positively defined.

Then one can find a projective subspace of dimension $(n-k-1)$
contained in the domain $U$ and a projective subspace of dimension
$k$ in the complement $\Bbb RP^n\setminus\overline U$.}
\endproclaim

{\tc Our Main Hypothesis and the very notion of
$L$-convex-concavity were invented during  an attempt to prove or
disprove  the Arnold hypothesis. We didn't succeed to prove it in
full generality. However, we obtained several results in this
direction.

We proved Arnold hypothesis for hypersurfaces satisfying the
following additional assumption: there exists a non-degenerate
quadratic cone  $K$ and a  hyperplane $\pi\subset \Bbb RP^n$ not
passing through the vertex of the cone, such that, first, the
hypersurface and the cone $K$ have the same intersection with the
hyperplane $\pi$, and, second, at each point of this intersection
the tangent planes to the hypersurface and to the cone coincide
(paper in preparation).

There is an {\tci affine} version of the Arnold hypothesis: one
should change  $\Bbb RP^n$ to $\Bbb R^n$ in its formulation (and
ask if there exist affine subspaces of dimensions $k$ and
$(n-k-1)$ in $U$ and $\Bbb R^n\setminus \overline U$
respectively). Our second result is an explicit construction of a
counterexample to this affine version of Arnold conjecture (paper
in preparation). The main role in this construction is played by
{\tci affine} convex-concave sets.

Here is the definition of the class of $(L)$-convex-concave
subsets of $\Bbb R^n$. Fix a class $(L)$ of $(k+1)$-dimensional
affine subspaces of  $\Bbb R^n$ parallel to $L$. Its elements are
parameterized by points of the quotient space  $\Bbb R^n/N$, where
$N$ is the (only) linear subspace of this class. A set $A$ is
called {\tci affine $(L-)$convex-concave} if
\item{1)} any section $A\cap N$  of $A$ by a subspace $N\in (L)$ is
convex and
\item{2)} the section $A\cap N_a$ depends  concavely on
the parameter $a\in \Bbb R^n/N$.

The last condition means that for any segment $a_t=ta+(1-t)b$,
$0\leq t\leq1$, in the parameter space $\Bbb R^n/N$ the section
$A\cap N_{a_t}$ is contained inside the linear combination (in the
Minkowski sense)  $t(A\cap N_a)+(1-t)(A\cap N_b)$ of the sections
$A\cap N_a$ and $A\cap N_b$. Any projective $L$-convex-concave set
is affine $(L)$-convex-concave in any affine chart not containing
$L$ with respect to  the class $(L)$ of  $(\dim L+1)$-dimensional
affine subspaces whose closures in $\Bbb RP^n$ contain  $L$.

For a class  $(L)$ of parallel planes in $\Bbb R^3$ we constructed
a $(L)$-convex-concave set $A\subset\Bbb R^3$ not containing lines
with smooth and everywhere non-degenerate boundary. However, all
our attempts to modify the example in such a way that its closure
$\overline{A}\subset\Bbb RP^3$ will be $L$-convex-concave failed.
Finally we proved that this is impossible: the Main Hypothesis is
true for $\Bbb R^3$  and any $L$-convex-concave set with $\dim
L=1$. This is the only case of the Main Hypothesis we were able to
prove (except trivially true cases of $\dim L=0$ and $\dim L=n-1$
in projective space $\Bbb RP^n$ of any dimension $n$).}

\subhead{\tcb The Main Hypothesis in the three-dimensional case}
\endsubhead
{\tc Our proof of the Main hypothesis in three-dimensional case is
quite lengthy. In this paper we construct an $L$-duality needed
for the fourth step of the proof (see below). The third step of
the proof requires a cumbersome combinatorics and will be given in
a separate paper.

We will give a sketch of this proof and will clarify the role of
$L$-duality.

{\tcr Sketch of the proof.} Any line lying inside a
$L$-convex-concave  set $A\subset \Bbb RP^3$ intersects all convex
sections $A\cap N$ of $A$ by planes $N$ containing the line $L$.
Vice versa, any line intersecting all these sections lies in $A$.
{\tcr The first step } of the proof is an application of a Helly
theorem \cite{He1, He2}. Consider a four-dimensional affine space
of all lines in $\Bbb RP^3$ not intersecting $L$, and convex
subsets $U_N$ of this space consisting of all line intersecting
the section  $A\cap N$. Applying the Helly theorem to the family
$U_N$, we conclude that if for any five sections $A\cap N_i$,
$i=1,\dots,5$, one can find a line intersecting all of them, then
there is a line intersecting all sections.

For any four section one can prove existence of a line
intersecting all of them. {\tcr The second step} of the proof
consists of the proof of this claim (in any dimension). }

\proclaim{\tcr Proposition 1 {\tc (about four sections)}}{\tci Let
$A$ be a $L$-convex-concave subset of $\Bbb RP^n$, and let $\dim
L=n-2$. Then for any four sections $A\cap N_i$ of the set $A$ by
hyperplanes $N_i$, $i=1,\dots,4$, $N_i\supset L$, one can find a
line intersecting all of them.}
\endproclaim

{\tc The proof uses a theorem due to Browder \cite{Br}. This
theorem is a version of a Brawer fixed point theorem claiming
existence of a fixed point of a continuous map of a closed
$n$-dimensional ball into itself. The Browder theorem deals with
set-valued upper semi-continuous maps of a convex set $B^n$ into
the set of all its closed convex subsets of $B^n$. The Browder
theorem claims that there is a point  $a\in B^n$ such that $a\in
f(a)$.

Here is how we use it. From the $L$-convex-concavity property of
the set $A\subset\Bbb RP^n$ with $\codim L=2$, one can easily
deduce that for any three sections $A_i=A\cap N_i$, $i=1,2,3$, and
any point $a_1\in A_1$ there is a line passing through  $a_1$ and
intersecting both $A_2$ and $A_3$. For four sections $A_i=A\cap
N_i$, $i=1,\dots,4$, and a point $a_1 \in A_1$ consider all pairs
of lines $l_1$ and $l_2$ such that
\item{1)} the line $l_1$ passes through $a_1$ and intersects $A_2$ and $A_3$,
\item{2)}  the line $l_2$ passes through the point of intersection of
$l_1$ and $A_3$, intersects  $A_4$ and intersects  $A_1$ at point
$a_1'$.

Consider a set-valued mapping $f$ of the section $A_1$ to the set
of all its subsets mapping the point  $a_1$ to the set of all
points $a_1'$ obtainable in this way. We prove that $f$ satisfies
conditions of the Browder theorem. Therefore  there exists a point
$a_1\in A_1$ such that $a_1\in f(a_1)$. It means that there is a
line $l_1$ passing through this point and coinciding  with the
corresponding line $l_2$. Therefore this line intersects the
sections $A_2$, $A_3$, $A_4$ and the second step of the proof ends
here.

Proof of the existence of a line intersecting (fixed from now on)
sections $A\cap N_i$, $i=1,\dots,5$, is quite complicated and goes
as follows. Choose an affine chart containing all five sections
and not containing the line $L$. Fix a Euclidean metric in this
chart.

Define a {\tci distance } from a line $l$  to the collection of
sections $A\cap N_i$, $i=1,\dots,5$, as the maximum of distances
from the point $a_i=l\cap N_i$ to the section $A\cap N_i$,
$i=1,\dots,5$. A line $l$ is a  {\tci Chebyshev} line if the
distance from $l$ to the sections $A\cap N_i$, $i=1,\dots,5$, is
the minimal one. We prove that for the Chebyshev line these
distances are all equal. With the Chebyshev line $l$ one can
associate five half-planes $p^{+}_i\subset N_i$. These half-planes
are supporting to the sections $A\cap N_i$ at the points $b_i\in
A\cap N_i$, the closest to $a_i$ points of the section $A\cap
N_i$. We have to prove that the distance from $L$ to the sections
is equal to zero, i.e. that $a_i=b_i$.

To prove it is enough to find a line $l'$ intersecting all
half-planes $p _{i}^{+}$, $i=1,\dots,5$. Indeed, if $a_i\neq b_i$
then, moving slightly the line $l$ into  the direction of the line
$l'$, one can decrease the distance from the line $l$ to the
sections  $A\cap N_i$, $i=1,\dots, 5$, which is impossible. So, it
is enough to prove that there exists a line $l$ intersecting the
five support half-planes $p^{+}_i \subset N_i$, $i=1,\dots,5$.

We will call the configuration of the five half-planes
$p^{+}_i\subset N_i$, $i=1,\dots,5$, {\tci non-degenerate} if
their boundaries intersect the line  $L$ in five different points.
Otherwise, i.e. if they intersect $L$ in less than five points, we
will  call the configuration  {\tci degenerate}. We prove the
existence of the line $l'$ separately for non-degenerate (Step 3)
and  degenerate (Step 4) cases.

Detailed proof of the third step is given in our paper  {\rm ``A
convex-concave domain in $\Bbb RP^3$ contains a line''} (in
preparation).

Here is a brief sketch of this {\tcr third step}. The proof of an
existence of a line intersecting all five half-planes
$p^{+}_i\subset N_i$ of a non-degenerate configuration is based on
a detailed analysis of combinatorial properties of each possible
configuration. It turns out that there are essentially only six
possible combinatorial types. For different combinatorial types of
configurations the proofs differ, though share the same spirit.

Here is a rough description of the most common scheme. Instead of
half-planes $p^{+}_i\subset N_i$, $i=1,\dots,5$, consider extended
half-planes $p_i$ such that
\item{1)} $p^{+}_i\subset p_i\subset N_i$;
\item{2)} boundaries of the half-planes $p_i$ intersect the Chebyshev line
and
\item{3)} intersections of the boundaries of $p_i$ and $p^{+}_i$ with the line $L$ coincide.

It is enough to prove that there exists a line intersecting all
extended half-planes $p_1\i$, $i=1,\dots,5$, and  at least one of
them at an interior point. Take planes $\pi _i$ containing the
Chebyshev line $l$ and boundaries of half-planes  $p_i$,
$i=1,\dots,5$. Each half-plane $p_j$ is divided by planes $\pi_i$
into five sectors. The minimizing property of the Chebyshev line
$l$ implies that some particular sectors necessarily intersect the
convex-concave set $A$.

Using combinatorial  properties of the configuration, we choose
four half-planes and a particular sector on one of them
intersecting the set $A$. Applying the Browder theorem (as on the
step 2), we prove existence of a line intersecting the four
sections {\it in some prescribed  sectors} of the corresponding
half-planes. From the combinatorial properties of the
configuration follows that the constructed line intersects the
fifth half-plane, q.e.d.

In the present paper we prove, among other results, the claim of
the  {\tcr fourth step}, i.e. existence of a line intersecting all
five half-planes $p^{+}_i\subset N_i$ of a degenerate
configuration. The proof goes as follows. All hyperplanes
$N\subset\Bbb RP^n$ containing a fixed subspace $L$ of codimension
2, can be parameterized by points of a projective line $\Bbb
RP^n/L$, so have a natural cyclic order. We say that a
$L$-convex-concave set $A$ with  $\dim L=n-2$ is {\tci linear
between cyclically ordered sections} $A_i=A\cap N_i$ if the
intersection $A_{ij}$ of the set $A$ with  a half-space of the
projective space bounded by two adjacent hyperplanes  $N_i$ and
$N_{i+1}$ coincides with  a convex hull of the sections $A_i=A\cap
N_i$ and $A_{i+1}=A\cap N_{i+1}$ (the convex hull is taken in any
affine chart $\Bbb RP^n\setminus N_j$, $j\neq i$, $j\neq i+1$, and
doesn't depend on the choice of the chart.)}

\proclaim{\tcr Proposition 2}{\tci Let $A$ be a $L$-convex-concave
subset of $\Bbb RP^n$, and $\dim L=n-2$. Suppose that   there
exist four sections of the set $A$ such that $A$ is linear between
these sections. Then the set $A$ contains a line.}
\endproclaim

{\tc This is a reformulation of the Proposition 1.

We prove the following, dual to the Proposition 2, claim.}

\proclaim{\tcr Proposition 3 {\tc (about sets with octagonal
sections)}}{\tci Let $D\subset\Bbb RP^n$ be a $L$-convex-concave
set, and $\dim L=1$. Suppose that any section $D\cap N$ of $D$ by
any  two-dimensional plane $N$ containing the line $L$, is an
octagon whose sides lie on lines intersecting the line $L$ in four
fixed (i.e. not depending on $N$) points. In other words, each
octagon has four pairs of "parallel" sides intersecting $L$ in a
fixed point. Then there exists an $(n-2)$-dimensional projective
subspace intersecting all planar sections $D\cap N$, $L\subset N$,
of the set $D$.}
 \endproclaim

{\tc  In fact, the main goal of this paper is to give a definition
of  an $L$-duality with respect to which the two propositions
above are dual, and to establish general properties of this
duality required for reduction of  the Proposition 3 to the
Proposition 2.

Let's return to the Step  4 of the proof. In  degenerate cases the
boundaries of the five half-planes  $p^{+}_i$, $i=1,\dots,5$
intersect the line $L$ in at most four points. Assume that their
number is exactly four and  denote them by $Q_1,Q_2,Q_3,Q_4$.
Perform now the following surgery of the set $A$. Replace each
convex section $A\cap N$ of the set $A$, $L\subset N$, by a
circumscribed octagon whose four pairs of parallel sides intersect
the line $L$ at the points  $Q_1,\dots,Q_4$. In \S 6 we prove that
application of this surgery to  a $L$-convex-concave set $A$
results in a $L$-convex-concave set $D$. The set $D$ satisfies
conditions of the Proposition 3, so there exists a line
intersecting all octagonal sections of the set $D$. This line
intersects all half-planes  $p^{+}_i$, $i=1,\dots,5$, and the
proof of the main hypothesis in three-dimensional case is
finished.}

\subhead{\tcb $L$-duality and plan of the paper}
\endsubhead
{\tc There are several well-known types of duality, e.g. a usual
projective duality or a duality between convex subsets of $\Bbb
R^n$ containing the origin and convex subsets of the dual space.
Different types of duality are useful for different  purposes.
Here we will construct a $L$-duality mapping a $L$-convex-concave
subset $A$ of a projective space $\Bbb RP^n$ to a set $A^\bot _L$
in the dual projective space  $(\Bbb RP^n)^*$. The set  $A^\bot
_L$ turns out to be $L^*$-convex-concave , where $L^*\subset(\Bbb
RP^n)^*$ is a subspace dual to $L$. The main duality property
holds for $L$-duality: $A=(A^\bot _L)^\bot _{L^*}$. The Main
Hypotheses for a set $A\subset \Bbb RP^n$ and for its dual $A^\bot
_L\subset (\Bbb RP^n)^*$ turn out to be equivalent: if the set
$A^\bot _L$ contains a projective subspace $M^*$ such that $\dim
M^*+\dim L^*=n-1$, then the set $A$ contains the dual subspace $M$
such that $\dim M+\dim L=n-1$. This is why $L$-duality is useful
for us: the problem for the $L$-dual set may be easier than for
the initial set. This is how the $L$-duality is used in the Step 4
of the proof of the Main Hypothesis in three-dimensional case.

In this paper we give a detailed description of  the $L$-duality.
Its meaning is easy to understand if the $L$-convex-concave set
$A$ is a domain with a smooth boundary. Assume that the boundary
$B$ of $A$ is strictly convex-concave, i.e. that its second
quadratic form is nondegenerate at each point. Consider a
hypersurface $B^*$ in the dual projective space $(\Bbb RP^n)^*$
projectively dual (in the classical sense) to $B$. The smooth
hypersurface $B^*$ divides $(\Bbb RP^n)^*$ into two parts. The
subspace  $L^*$ dual to $L$ doesn't intersect hypersurface $B^*$,
so exactly one of the connected component of $(\Bbb
RP^n)^*\setminus B^*$ does not contain $L^*$. The $L$-dual of the
set $A$ coincides with the closure of this component.

This definition does not work for sets whose boundary is not
smooth and strictly convex-concave. However, we are forced to deal
with such sets (in particular with sets whose sections are closed
convex polygons and whose complements to projections are open
convex polygons). Therefore we have to give a different, more
suitable to our settings definition. An example of how one can
define such a thing  is the classical definition of  dual convex
sets. We follow closely this example.

Here is the plan of the paper. First, in \S1, we give a definition
of projective separability, mimicking the standard definition of
separability for  affine spaces. All statements formulated in this
paragraph are immediate, so we omit the proofs. In \S2 we discuss
the notion of projective duality, the notion mimicking the
classical definition of duality for containing the origin convex
subsets of  linear spaces. Here all statements are also very
simple, but for the sake of completeness we give their proofs and
explain why all of them are parallel to the classical ones.

After that, in \S3, we define $L$-duality and prove its basic
properties (using already defined projective separability and
projective duality). At the end of \S3 we discuss semi-algebraic
$L$-convex-concave sets and a relation between the $L$-duality and
integration by Euler characteristics. The results of \S5 and \S6
will be used in the Step 4 of the proof of the Main Hypothesis in
the three-dimensional case. From the results of \S4 follows, in
particular, the proposition about convex-concave sets with
octagonal sections (the Proposition~3 above). In \S6 we describe,
in particular, the surgery allowing to circumscribe convex
octagons around planar convex sections.}


\head{\tcb \S1. Projective and affine separability}
\endhead

{\tc We recall the terminology related to the notion of
separability in projective and affine spaces.}

\subhead{\tcb Projective case}
\endsubhead
{\tc We say that a subset $A\subset \Bbb RP^n$ is  {\tci
projectively separable} if any point of its complement lies on a
hyperplane not intersecting the set  $A$. }

\proclaim{\tcr Proposition}{\tci Complement to a projectively
separable set $A$ coincides with  a union of all hyperplanes not
intersecting the set $A$. Vice versa, complement to any union of
hyperplanes has property of projective separability.}
\endproclaim

{\tc This proposition can be reformulated:}

\proclaim{\tcr Proposition}{\tci Any subset of  projective space
defined by a system of linear homogeneous inequalities
$L_\alpha\neq 0$, where $\alpha$ belongs to some set of indexes
and  $L_\alpha $ is a homogeneous polynomial   of degree one, is
projectively separable. Vice versa, any projectively separable set
can be defined in this way.}
\endproclaim

{\tc  We define a   {\tci projective separability hull of the set}
$A$ as the smallest projectively separable set containing the set
$A$.}

\proclaim{\tcr Proposition}{\tci The projective separability hull
of a set  $A$ is exactly the complement to a union of all
hyperplanes  in  $\Bbb RP^n$ not intersecting the set $A$. In
other words, a point lies in the projective separability hull of
the set $A$ if and only if any hyperplane containing this point
intersects the set $A$.}
\endproclaim

\subhead{\tcb Affine case}
\endsubhead
{\tc Recall the well known notion of separability in  the affine
case. Namely, a subset $A$ of an affine space is {\tci affinely
separable} if any point of the complement to the set $A$ belongs
to a closed half-space not intersecting the set $A$. Evidently,
any affinely separable set is convex and connected.}

\proclaim{\tcr Proposition}{\tci The complement to an
affinely separable set $A$ coincides with  a union of closed
half-spaces not intersecting the set $A$. Vice versa, a complement
to any union of closed half-spaces is affinely separable.}
\endproclaim

{\tc This property can be reformulated.}

\proclaim{\tcr Proposition}{\tci Any subset of an affine space
defined by a system of linear inequalities $\{L_\alpha (x)< 0\}$,
where $\alpha$ belongs to some set of indices and  $L_\alpha $ is
a polynomial   of degree at most one, is affinely separable. Vice
versa, any affinely separable set  can be defined in this way.}
\endproclaim

{\tc  We define an   {\tci affine separability hull} of a set $A$
as the smallest set containing the set $A$ and having the property
of affine separability.}

\proclaim{\tcr Proposition}{\tci Affine separability hull of a set
$A$ is equal to a complement to a union of all closed not
intersecting the set $A$ half-spaces of the affine space. In other
words, a point lies in the affine separability hull of the set $A$
if and only if any closed half-space containing this point also
intersects the set $A$.}
\endproclaim

\subhead{\tcb Convex subsets of projective spaces and
separability}
\endsubhead
{\tc  Projective and affine separability are closely connected.}

\proclaim{\tcr Proposition}{\tci Let $L$ be a
 hyperplane in a projective space $\Bbb RP^n$ and     $U=\Bbb
RP^n\setminus L$ be a corresponding affine chart.

1. Any affinely separable subset of the affine chart $U$  (so, in
particular, connected and convex in $U$),  is also projectively
separable as a subset of a projective space.

2. Any connected projectively separable subset of the affine chart
$U$ is  also affinely separable as a subset of an affine space
$U$.}
\endproclaim

{\tc A connected projectively separable subset of a projective
space  not intersecting at least one hyperplane will be called a
{\tci separable convex} subset of the projective space. (There is
exactly one projective separable subset of  projective space
intersecting all hyperplanes, namely the projective space
itself.)}

\remark{\tcr Remark}{\tc We defined above a notion of a (not
necessarily projectively separable) convex subset of a projective
space: a nonempty  subset $A$ of a projective space $\Bbb RP^n$ is
called {\tci convex} if, first, there is a hyperplane
$L\subset\Bbb RP^n$ not intersecting the set $A$ and, second, any
two points of the set $A$ can be joined by a  segment lying in
$A$. We will not need convex non-separable sets.}
\endremark

\head{\tcb \S2. Projective and linear  duality }
\endhead

{\tc We construct here a variant of a projective duality. To a
subset $A$ of a projective space $\Bbb RP^n$ corresponds in virtue
of this duality a subset $A^*_p$ of the dual projective space
$(\Bbb RP^n)^*$. This duality is completely different from the
usual projective duality and is similar to a linear duality used
in convex analysis. For the sake of completeness we describe here
this parallelism as well.}

\subhead{\tcb Projective duality}
\endsubhead
{\tc Projective space $\Bbb RP^n$ is obtained as a factor of a
linear space $\Bbb R^{n+1}\setminus 0$ by a proportionality
relation. The dual projective space, by definition, is a factor of
the  set of all nonzero covectors  $\alpha\in (\Bbb
R^{n+1})^*\setminus 0$ by a proportionality relation.

There is  a one-to-one correspondence between hyperplanes in the
space  and points of the dual space. More general, to any subspace
$L\subset \Bbb RP^n$ corresponds a dual subspace  $L^*\subset
(\Bbb RP^n)^*$ of all hyperplanes containing  $L$, and
the duality property $(L^*)^*=L$ holds.

For any set  $A\in\Bbb RP^n$ we define its dual  set
$A^*_p\subset(\Bbb RP^n)^*$  to be a set of all hyperplanes in
$\Bbb RP^n$ not intersecting the set $A$. (The symbol  $A^*$
denotes the dual space, so we introduce the new notation
$A^*_p$.)}

\proclaim{\tcr Proposition}{\tci 1. If $A$ is non-empty, then the
set  $A^*_p$ is contained in some affine chart of the dual space.

2. The set $A^*_p$ is projectively separable.}
\endproclaim

\demo{\tcr Proof}{\tc 1. The set $A$ is nonempty, so contains some
point $b$. A hyperplane $b^*\in (\Bbb RP^n)^*$ corresponding to
the point $b$, doesn't intersect the set $A^*_p$. Therefore the
set $A^*_p$ is contained in the  affine chart $(\Bbb
RP^n)^*\setminus b^* $.

2. If a   hyperplane $L\subset\Bbb RP^n$, considered as a point in
the space $(\Bbb RP^n)^*$, is not contained in the set $A^*_p$,
then, by definition, the  hyperplane $L$ intersects the set $A$.
Let $b \in A\cap L$. The hyperplane $b^*$ dual to the point $b$
doesn't intersect the set $A^*_p$. So this hyperplane separates
the point corresponding to the hyperplane $L$ from the set
$A^*_p$.}
\enddemo

{\tc The following theorem gives a full description of the set
$(A^*_p)^*_p$.}

\proclaim{\tcr Theorem}{\tci For any  set $A\subset \Bbb RP^n$ the
corresponding set $(A^*_p)^*_p$ consists of all points $a$ such
that  any hyperplane containing  $a$ intersects the set $A$. In
other words, the set $(A^*_p)^*_p$ coincides with the projective
separability hull of the set $A$.}
\endproclaim

\demo{\tcr Proof}{\tc The point  $a$ belongs to $(A^*_p)^*_p$ if
and only if the corresponding hyperplane $a^*\subset(\Bbb RP^n)^*$
doesn't intersect the set $A^*_p$. To any point $p$ in $(\Bbb
RP^n)^*$ of this hyperplane corresponds a hyperplane
$p^*\subset\Bbb RP^n$ containing the point $a$. The point
$p\in(\Bbb RP^n)^*$ doesn't belong to $A^*_p$ if and only if the
hyperplane $p^*\subset\Bbb RP^n$ intersects the  set $A$. So the
condition that all points of the hyperplane $a\subset(\Bbb
RP^n)^*$ does not belong to $A^*_p$, means that all hyperplanes in
$\Bbb RP^n$ containing the point $a$, intersect the set $A$.}
\enddemo

\proclaim{\tcr Corollary}{\tci The duality property
$(A^*_p)^*_p=A$ holds for all projectively separable subsets of a
projective space, and only for them.}
\endproclaim

\subhead{\tcb Linear duality}
\endsubhead
{\tc The property of affine separability differs from the property
of projective separability: we use closed half-spaces in the
affine definition and hyperplanes in the projective definition.
One can do the same with the duality theory developed above and
define the set $A^*_a$ corresponding to a subset of an affine
space as a set of all closed half-spaces not intersecting the set
$A$. This definition is not very convenient because the set of all
closed half-spaces doesn't have a structure of an affine space.
Moreover, this set is topologically different from affine space:
it is homeomorphic to the sphere $S^{n}$ with two removed points
(one point corresponding to an empty set and another to the whole
space). One can avoid this difficulty by considering instead a set
of all closed half-spaces not containing some fixed point with one
added element (this element  corresponds to an empty set regarded
as a half-space on an infinite distance from the fixed point).
This set has a natural structure of an affine space. Namely,
taking the fixed point as the origin and denoting the resulting
linear space by $\Bbb R^n$, one can parameterize the set described
above by $(\Bbb R^n)^*$: to any nonzero $\alpha\in(\Bbb R^n)^*$
corresponds a closed half-space defined by inequality $\langle
\alpha,x \rangle \geq 1$. To $\alpha =0$ corresponds an empty set
(defined by the same inequality  $\langle \alpha,x \rangle \geq
1$).

It is more convenient to consider only sets containing some fixed
point when talking about affine duality. Taking this point as the
origin, we get the well-known theory of affine duality, which is
parallel to the theory of projective duality. Here are its main
points.

To any subset $A$ of a linear space $\Bbb R^n$ corresponds a
subset $A^*_l$ of a dual space $(\Bbb R^n)^*$ consisting of all
 $\alpha \in (\Bbb R^n)^*$ such that the inequality $\langle \alpha,x
\rangle <1$ holds for all $x\in A$.}

\proclaim{\tcr Proposition}{\tci For any  set $A\subset \Bbb R^n$
containing the origin the corresponding dual  set $A^*_l$ in the
dual space has the property of affine separability. In particular,
it is convex.}
\endproclaim

\proclaim{\tcr Proposition}{\tci  For any  set $A\subset \Bbb R^n$
containing the origin the  set $(A^*_l)^*_l$ consists of all
points $a\in\Bbb R^n$ with the following property: any closed
half-space containing  $a$ intersects the set $A$. In other words,
the  set $(A^*_l)^*_l$ is equal to the affine separability hull of
the  set $A$.}
\endproclaim

\proclaim{\tcr Corollary}{\tci The duality property
$(A^*_p)^*_p=A$ holds for all containing the origin convex sets
with the property of affine separability, and only for them.}
\endproclaim

\head{\tcb \S3. $L$-duality}
\endhead

{\tc Here we construct a $L$-duality. A subset $A$ of a projective
space $\Bbb RP^n$ disjoint from some subspace $L$, will be
$L$-dual to a subset $A^\bot _L$ of a dual projective space $(\Bbb
RP^n)^*$ disjoint from the subspace $L^*$.

Any subset $C$ in the projective space $(\Bbb RP^n)^*$ can be
considered as a subset of a set of all hyperplanes in the
projective space $\Bbb RP^n$. We will also denote it by $C$.

Let $L$ be some projective subspace of   $\Bbb RP^n$, and $A$ be
any set not intersecting $L$. For a hyperplane $\pi$ not
containing the subspace $L$, denote by $L_\pi$ the subspace $L\cap
\pi$. Consider a factor-space $(\Bbb RP^n)/L_\pi$. The image
$\pi_L$ of a hyperplane $\pi$ is a hyperplane in the factor-space
$(\Bbb RP^n)/L_\pi$.}

\definition{\tcr Definition}{\tc We say that the   hyperplane $\pi$
{\tci belongs to the $L$-dual set} $A^\bot_L$ if $\pi$ doesn't
contain $L$ and the hyperplane  $\pi_L$ is contained in the
projection of the set $A$ on the factor-space $(\Bbb
RP^n)/L_\pi$.}
\enddefinition

{\tc In other words, a hyperplane $\pi$ belongs to the set $A^\bot
_L$ if projection of $\pi$ from the center $L_\pi$ belongs to
$B^*_p$, where $B$ is the complement to the projection of the set
$A$ on the space $\Bbb RP^n/L_\pi$.

Here is another description of the  set $A^\bot_L$. The complement
$\Bbb RP^n\setminus L$ to the subspace $L$ is fibered by  spaces
$N\supset L$ of dimension $\dim N=\dim L+1$. A hyperplane $\pi$
belongs to $A^\bot_L$, if and only if   for any fiber $N$ its
intersection with the set $A\cap\pi$ is non-empty, $N\cap A\cap
\pi \neq \emptyset $.} In other words, $\pi\in A^\bot_L$ if and
only if $\pi$ intersects any section of $A$ by any $(\dim
L+1)$-dimensional space containing $L$.

\example{\tcr Example}{\tc Let $L$ be a hyperplane, and $A$ be a
set disjoint from $L$, $A\cap L=\emptyset$. Then $A^\bot_L$ is a
union of all hyperplanes intersecting the set $A$. In other words,
the set $A^\bot _L$ is a complement to the set $A^*_p$. Indeed, in
this case the only space $N$ containing $L$ is the projective
space $\Bbb RP^n$ itself. Note that in this case the set $L$-dual
to $A$ doesn't depend on the choice of a hyperplane $L$ (as long
as $L$ doesn't intersect the set $A$).}
\endexample

\proclaim{\tcr Proposition}{\tci If $A\subset B$ and $B\cap
L=\emptyset$, then $(A^\bot _L)\subset (B^\bot _L)$.}
\endproclaim

\demo{\tcr Proof}{\tc If a hyperplane intersects all sections
$A\cap N$, then it intersects all sections $B\cap N$.}
\enddemo

\proclaim{\tcr Proposition}{\tci Let $M$ be a projective subspace
in  $\Bbb RP^n$  not intersecting $L$ and of a maximal possible
dimension, i.e. $\dim M =\dim L^*= n-\dim L-1$. Then $M^\bot_L =
M^*$.}
\endproclaim

\demo{\tcr Proof}{\tc Any section of $M$ by  $(\dim
L+1)$-dimensional space containing $L$ is just a point, and any
point of $M$  is a section of $M$ by such a space. By definition
of $M^\bot_L$, a hyperplane $\pi$ belongs to  $M^\bot_L$ if and
only if it intersects any such section, i.e. contains any point of
$M$. This is exactly the definition of $M^*$. }
\enddemo

{\tc Let $L^*\subset (\Bbb RP^n)^*$ be the space dual to $L$. What
can be said about: a) sections of the set $A^\bot_L$ by $(\dim L^*
+1)$-dimensional  spaces $N\supset L^*$; b) projections of the set
$A^\bot _L$ from a $(\dim L^*-1)$-dimensional subspace $T$ of a
space $L^*$? We give below answers to these questions.}

\subhead{\tcb Sections of the  $L$-dual set}
\endsubhead
{\tc Recall first a duality between sections and projections. Let
$N$ be a projective subspace in the  space $(\Bbb RP^n)^*$.
Consider a  dual to $N$ subspace $N^*\subset\Bbb RP^n$. We will
need later an isomorphism and a projection described below.

There is a natural isomorphism between a space dual to the
quotient space  $\Bbb RP^n/N^*$ and the space $N$. This
isomorphism is a projectivization of a natural isomorphism between
a space dual to a factor-space and a subspace of a dual space dual
to the
 kernel of the factorization. Each hyperplane containing the
space $N^*$, projects to a hyperplane in $\Bbb RP^n/N^*$. (If a
hyperplane doesn't contain the space $N^*$, then its projection is
the whole space $\Bbb RP^n/N^*$.)

Using this isomorphism one can describe a section  of the  set
$C\subset (\Bbb RP^n)^*$  by the space $N$ in terms of the space
$\Bbb RP^n$. Consider a subset $C_{N^*}\subset(\Bbb RP^n)^*$ of
the set of hyperplanes $C$ consisting of all hyperplanes
containing $N^*$ (this is equivalent to $C_{N^*}=C\cap N$). Each
hyperplane from $C_{N^*}$ projects to a hyperplane in the
factor-space $\Bbb RP^n/N^*$. But the space $(\Bbb RP^n/N^*)^*$ is
identified with the space $N$. After projection and identifying we
get the required section $C\cap N$ from the set $C_{N^*}$.}

\proclaim{\tcr Theorem 1}{\tci Let $A$ be a subset of $\Bbb RP^n$
not intersecting $L$, and $N$ be any subspace of $(\Bbb RP^n)^*$,
containing $L^*$ as a hyperplane (i.e. $\dim N=\dim L+1$ and
$N\supset L^*$). Then the section $A^\bot _L\cap N$ is equal to
$B^*_p$, where $B \subset (\Bbb RP^n/N^*)$ is a complement to the
projection  of the set $A$ on the space $(\Bbb RP^n)/N^*$.}
\endproclaim

\demo{\tcr Proof}{\tc This  Theorem follows from the description
above of  sections  of subsets of $(\Bbb RP^n)^*$. Consider the
set of hyperplanes $C=A^\bot_L$. By the definition of the set
$A^\bot _L$, the set  $C_{N^*}$ consists of all hyperplanes
containing the projective space $N^*$, such that their projections
on $\Bbb RP^N/N^*$ after projection from $N^*$ are contained in
projection of the set $A$. In other words, their projections are
hyperplanes in $\Bbb RP^N/N^*$ not intersecting the complement to
the projection of the set $A$. Vice versa, any hyperplane not
intersecting this complement $B$ is, by definition of the set
$A^\bot _L$, a projection of some hyperplane belonging to the set
$C_{N^*}$. Therefore $A^\bot _L\cap N= B^*_p$.}
\enddemo

\subhead{\tcb Projections of $L$-dual sets}
\endsubhead
{\tc Recall a duality between projections and sections.

Denote by $Q$ a subspace in $\Bbb RP^n$ dual to the center of
projection $T\subset(\Bbb RP^n)^*$. There is a natural isomorphism
between the space $Q^*$,
 consisting of all hyperplanes of the space $Q$, and the factor-space
$(\Bbb RP^n)^*/T$. Namely, one should consider points of $(\Bbb
RP^n)^*/T$ as equivalency classes in the set of all hyperplanes in
the space $\Bbb RP^n$ not containing the space $Q$, of the
following equivalency relation: two hyperplanes are equivalent if
and only if their intersections with $Q$ coincide. This
intersection is the hyperplane in the space $Q$ corresponding
 to this equivalency class.

Projection of a  subset $C$ of $(\Bbb RP^n)^*$ from a center $T$
can be described in the following way. A set of hyperplanes $C$ in
$\Bbb RP^n$ defines some set of hyperplanes $C(Q)$ in the subspace
$Q=T^*$: a hyperplane $Q_1\subset Q$ belongs to the set $C(Q)$ if
and only if there exists a hyperplane belonging to the set $C$
intersecting $Q$ exactly by $Q_1$. Projection of the set $C$ from
the center $T$ is exactly the set $C(Q)$ of hyperplanes in $Q$
after identifying $Q^*$ and $(\Bbb RP^n)^*/T$.}

\proclaim{\tcr Theorem 2}{\tci Let $A$ be a set in
 $\Bbb RP^n$ not intersecting $L$, and $T$ be a
 hyperplane in the dual space $L^*\subset (\Bbb RP^n)^*$. Then the
projection of the set $A^\bot _L$ from the center $T$ can be
described as a set of all hyperplanes $p$ in space $Q=T^*\supset
L$ with the following property: there exists a  hyperplane $\pi
\subset A^\bot_L$ whose intersection with $Q$ is equal to $p$,
$p=\pi\cap Q$.}
\endproclaim

\demo{\tcr Proof}{\tc This Theorem follows from the description of
projections of subsets $C\subset(\Bbb RP^n)^*$ given  above.}
\enddemo

\definition{\tcr Definition}{\tc We say that a set $A$ is {\tci
coseparable} relative to $L$ if $A\cap L=\emptyset$ and for any
hyperplane $L_1 \subset L$ a complement to projection of the set
$A$ from the center $L_1$ has the property of affine separability
in space $(\Bbb RP^n)/L_1$.}
\enddefinition

\proclaim{\tcr Corollary}{\tci If, in addition to all conditions
of the Theorem~2, the set $A$ is coseparable relative to $L$, then
the complement to the projection of  the set $A^\bot_L$ from the
center $T$ is dual to the section  $A\cap T^* $ (i.e. equal to
$(A\cap T^*)^*_p$).}
\endproclaim

\subhead{\tcb Description of the set  $(A^\bot
_L)^\bot_{L^*}$}
\endsubhead
{\tc Let $A$ be a subset of $\Bbb RP^n$ not intersecting a
subspace $L$, and $L^*$ be a dual to $L$ subspace of $(\Bbb
RP^n)^*$. What can be said about a subset of  $\Bbb RP^n$
$L^*$-dual to the subset  $A^\bot _L$ of the space $(\Bbb
RP^n)^*$? From the theorems~1 and 2 we easily obtain the
description of this set $(A^\bot _L)^\bot_{L^*}$.}

\proclaim{\tcr Theorem 3}{\tci The set $(A^\bot _L)^\bot_{L^*}$
doesn't intersect $L$ and  consists of all points $a\in\Bbb RP^n$
satisfying the following condition: in the space  $L_a$ spanned by
$L$ and $a$, for any hyperplane $p$ in $L_a$ containing the point
$a$, $a\in p\subset L_a$ there is a hyperplane $\pi\subset \Bbb
RP^n$, $\pi\in A^\bot _L$, such that $p=\pi\cap L_a$.}
\endproclaim

\demo{\tcr Proof}{\tc  A section of the set $(A^\bot _L)^\bot
_{L^*}$ by the subspace $L_a$ can be described, according to the
theorem ~1 (applied to the subset $A^\bot _L$ of the space $(\Bbb
RP^n)^*$ and the subspace $L^*$ of this space), as the set of
hyperplanes in the factor-space $(\Bbb RP^n)^*/L^*_a$ not
intersecting a complement to the projection of the set $A^\bot _L$
on the space $(\Bbb RP^n)^*/L^*_a$.

So the point $a\in \Bbb RP^n$ lies in $(A^\bot_L)^\bot_{L^*}$ if
and only if a hyperplane in $(\Bbb RP^n)^*/L^*_a$, corresponding
to this point $a\in \Bbb RP^n$, $a\in L_a$, is contained in the
projection of the set $A^\bot_{L}$. This means that any hyperplane
$p$ of $L_a$, $P\subset L_a$, containing the point $a$, lies in
the projection of the set  $A^\bot_L$, if considered as  a point
of the space  $(\Bbb RP^n)^*/L^*_a$. This means, according to the
theorem~2, that for the hyperplane $p$ there exists a hyperplane
$\pi \in A^\bot_l$ such that $\pi \cap L_a=p$, q.e.d.}
\enddemo

{\tc Let's reformulate the Theorem~3. The point  $a$ belongs to
the set $(A^\bot _L)^\bot_{L^*}$ if the following two conditions
hold:}

\remark{\tcr Condition 1}{\tc The point $a$ in the space $L_a$,
spanned by $L$ and $a$, has the following property: any
 hyperplane $p\subset L_a$, containing the point
$a$, intersects the set $L_a\cap A$. In other words, the point $a$
lies in the set $((L_a\cap A)^*_p)^*_p$.}
\endremark

\remark{\tcr Condition 2}{\tc Projection  of the point $a$ from
any center $L_1\subset L$, where $L_1$ is a hyperplane in $L$, is
contained in some hyperplane in the space $(\Bbb RP^n)/L_1$
contained in the projection of the set $A$ on the space $(\Bbb
RP^n)/L_1$.}
\endremark

\proclaim{\tcr Theorem 4}{\tci The conditions 1 and 2 are
equivalent to the condition that the point $a$ belongs to the set $(A^\bot _L)^\bot_{L^*}$.}
\endproclaim

\demo{\tcr Proof}{\tc Indeed , according to the Theorem~3, if
$a\in (A^\bot _L)^\bot_{L^*}$, then any hyperplane $p$ in the
space $L_a$ containing the point $a$, is an intersection of $L_a$
and a hyperplane $\pi\in A^\bot_L$. This means that, first, the
hyperplane $p$ intersects $A$ and, second, that the projection of
the point $a$ from $L_1=L\cap \pi$ is containing in a hyperplane
in the factor-space $\Bbb RP^n/L_1$, which, in turn, is contained
in the projection of the set $A$. The first property is equivalent
to the Condition~1, and the second is equivalent to the
Condition~2.}
\enddemo

\proclaim{\tcr Corollary}{\tci Suppose that a  set $A$ doesn't
intersect the space $L$, and  intersection of $A$ with any
subspace $N$ containing $L$ as a hyperplane, is projectively
separable in projective space $N$. Then $(A^\bot
_L)^\bot_{L^*}\subset A$.}
\endproclaim

\demo{\tcr Proof}{\tc Indeed, the Condition~1 guarantees that for
any space $N$, containing $L$ as a hyperplane, the inclusion
$(A^\bot _L)^\bot_{L^*}\cap N\subset ((N\cap A)^*_p)^*_p$ holds.
But $((N\cap A)^*_p)^*_p =N\cap A$, since $N\cap A$ is
projectively separable. Therefore $(A^\bot _L)^\bot_{L^*}\subset
A$.}
\enddemo

\proclaim{\tcr Corollary}{\tci Suppose that the set $A$ is
coseparable relative to $L$. Then the intersection of the set
$(A_L^\bot)^\bot_{L^*}$ with any space $N$, containing $L$ as a
hyperplane, depends on the subset $A\cap N$ of the projective
space $N$ only and coincides with the set $((A\cap N)^*_p)^*_p$.
In particular, in this case $A\subseteq
(A_L^\bot)^\bot_{L^*}$.}
\endproclaim

\demo{\tcr Proof}{\tc If the set $A$ is coseparable relative to
$L$, then the Condition~2 holds for points satisfying to the
Condition~1. This is exactly what the Corollary claims. }
\enddemo

\subhead{\tcb Properties of $L$-coseparable and $L$-separable
sets}
\endsubhead
{\tc Let's sum up the facts about $L$-coseparable and
$L$-separable subsets of a projective space proved above.

Let a subset $A$ of a  projective space $\Bbb RP^n$ be coseparable
relative to a space $L$, and  suppose that any section of $A$ by a
space containing $L$ as a hyperplane, is projectively separable.

Then the set $A^\bot _L$  in the dual projective space $(\Bbb
RP^n)^*$ has the same properties relative to the dual space $L^*$.
Moreover, any section of $A^\bot_L$ by a subspace $N$ containing
$L^*$ as a hyperplane, is dual to the set $B$ (i.e. is equal to
$B^*_p$), where $B$ is a complement to the projection of the set
$A$ on $(\Bbb RP^n)/N^*$ from the center $N^*$. Projection of the
set $A^\bot _L$ from the center $T$, where $T$ is any hyperplane
in space $L^*$, is dual to the section of $A$ by $T^*$ (i.e. is
equal to $(A\cap T^*)^*_p$). Also, the duality relation $(A^\bot
_L)^\bot _L=A$ holds.

If the set $A^\bot _L$ contains a projective space $M^*$ of
dimension equal to the dimension of the  space $L$, then the set
$A$ contains its dual space $M$ of dimension equal to the
dimension of the space $L^*$.

 $L$-convex-concave sets  are
$L$-separable and $L$-coseparable, because closed sets and open
sets are both separable. Therefore for $L$-convex-concave all the
aforementioned properties hold. }

\subhead{\tcb Semialgebraic $L$-convex-concave  sets}
\endsubhead
{\tc Here we will use the integration by Euler characteristics,
introduced by O.~Viro (see [Vi]). We will denote
  Euler characteristics of a set $X$ by $\chi (X)$.}

\proclaim{\tcr Theorem}{\tci Let $A$ be a $L$-convex-concave
closed semialgebraic set in $\Bbb RP^n$, and let $\dim L=k$. Then
for any hyperplane $\pi\subset \Bbb RP^n$ the $\chi(A\cap \pi)$ is
equal to $\chi (\Bbb RP^{n-k-1})$ or to $\chi (\Bbb RP^{n-k-2})$.
In the first case the hyperplane $\pi$, considered as a  point of
$(\Bbb RP^n)^*$, belongs to the $L$-dual to $A$ set $A^\bot_L$. In
the second case the hyperplane $\pi$ doesn't belong to the set
$A^\bot _L$.}
\endproclaim

\demo{\tcr Proof}{\tc The complement to $L$ in $\Bbb RP^n$ is a
union of nonintersecting fibers, each fiber being a
$(k+1)$-dimensional space $N$ containing $L$. The set $A$ is
$L$-convex-concave, so its intersection with each fiber $N$ is
convex and closed. Therefore for each space $N$ the intersection
$A\cap N\cap \pi$ of the set $A\cap N$ with a hyperplane $\pi$
either is empty or is a closed convex set.

Suppose that the hyperplane $\pi$ doesn't contain the space $L$,
and denote by $L_\pi$ the space $L\cap \pi$. In the factor-space
$\Bbb RP^n/L_\pi$ we have a fixed point $\pi(L)$ (projection  of
the space $L$), a set $B$ (the complement to the projection of the
set $A$ from $L_{\pi}$), and a  hyperplane $\pi _L$ (projection of
the hyperplane $\pi$). To each point $a$ of the hyperplane $\pi_L$
in the factor-space corresponds a space $N(a)$ in $\Bbb RP^n$,
$N(a)\supset L$, whose projection is equal to the line passing
through $a$ and $\pi(L)$. The intersection $N(a)\cap A\cap \pi$ is
empty if $a$ belongs to the set $B$. Otherwise, the intersection
$N(a)\cap A\cap \pi$ is a closed convex set. The Euler
characteristics of the set $N(a)\cap A\cap \pi$ is equal to zero
in the first case, and is equal to one in the second case. Using
Fubini theorem for an integral by Euler characteristics for the
projection of the set $A\cap \pi$ on the factor-space  $\Bbb
RP^n/L_\pi$, we get
$$
\chi (A\cap \pi) =\chi
(\pi_L\setminus (\pi_L\cap B)).
 $$
So  $\chi (A\cap \pi)=\chi(\pi_L)=\chi (\Bbb RP^{n-k-1})$, if
$\pi_L\cap B=\emptyset$. Otherwise, i.e. if  $\pi_L\cap B\neq
\emptyset$, the $\chi (A\cap\pi)=\chi (\Bbb RP^{n-k-2})$. In the
first case $\pi_L\in A^\bot _L$ by definition, and in the second
case $\pi_L\notin A^\bot_L$. Therefore  the theorem is proved for
hyperplanes not containing the space $L$. If $L\subset\pi$, then
from similar considerations one can see that $\chi (\pi\cap
A)=\chi(\Bbb RP^{n-k-2})$, q.e.d.}
\enddemo

\proclaim{\tcr Corollary}{\tci For  semi-algebraic
$L$-convex-concave sets $A$ the $L$-dual set is defined
canonically (i.e. $A^\bot_L$ doesn't depend on the choice of the
space $L$, relative to which the set $A$ is $L$-convex-concave).}
\endproclaim

\remark{\tcr Remark}{\tc For the semialgebraic $L$-convex-concave
sets one can prove the duality relation
$$
(A^\bot _L)^\bot _{L^*}=A,
$$
using only this theorem and a Radon transform for the integral by
Euler characteristics, see [{\rm Vi}], and also  [{\rm PKh}].}
\endremark

\head{\tcb \S4. Duality between pointed convex sections of
convex-concave  sets and affine dependence of convex sections on
parameter}
\endhead

{\tc In  this section we define properties of pointedness (with
respect to a cone)  and of affine dependence on parameter (for
parameters belonging to some convex domain) of sections.

We begin with affine versions of these notions and then give
corresponding projective definitions. We prove that the property
of pointedness and the property of affine dependence on parameters
are dual.

\subhead{\tcb Pointedness of sections}
\endsubhead
{\tc We start with affine settings.  Let $K$ be a pointed
 (i.e. not containing  linear subspaces) closed
convex cone in a linear space $N$ with vertex at the origin.

We say that a set $A$ is {\tci pointed} with respect to  $K$, if
there is a point $a\in A$ such that the set $A$ lies entirely in a
translated cone $(K+a)$ with the vertex at the point $a$. This
point $a$ will be called a {\tci vertex} of the set $A$ relative
to the cone $K$. The vertex of the set $A$ relative to $K$ is
evidently uniquely defined.

\midinsert
\input epsf
\centerline{\hfill\epsfysize=0.13\vsize\epsffile{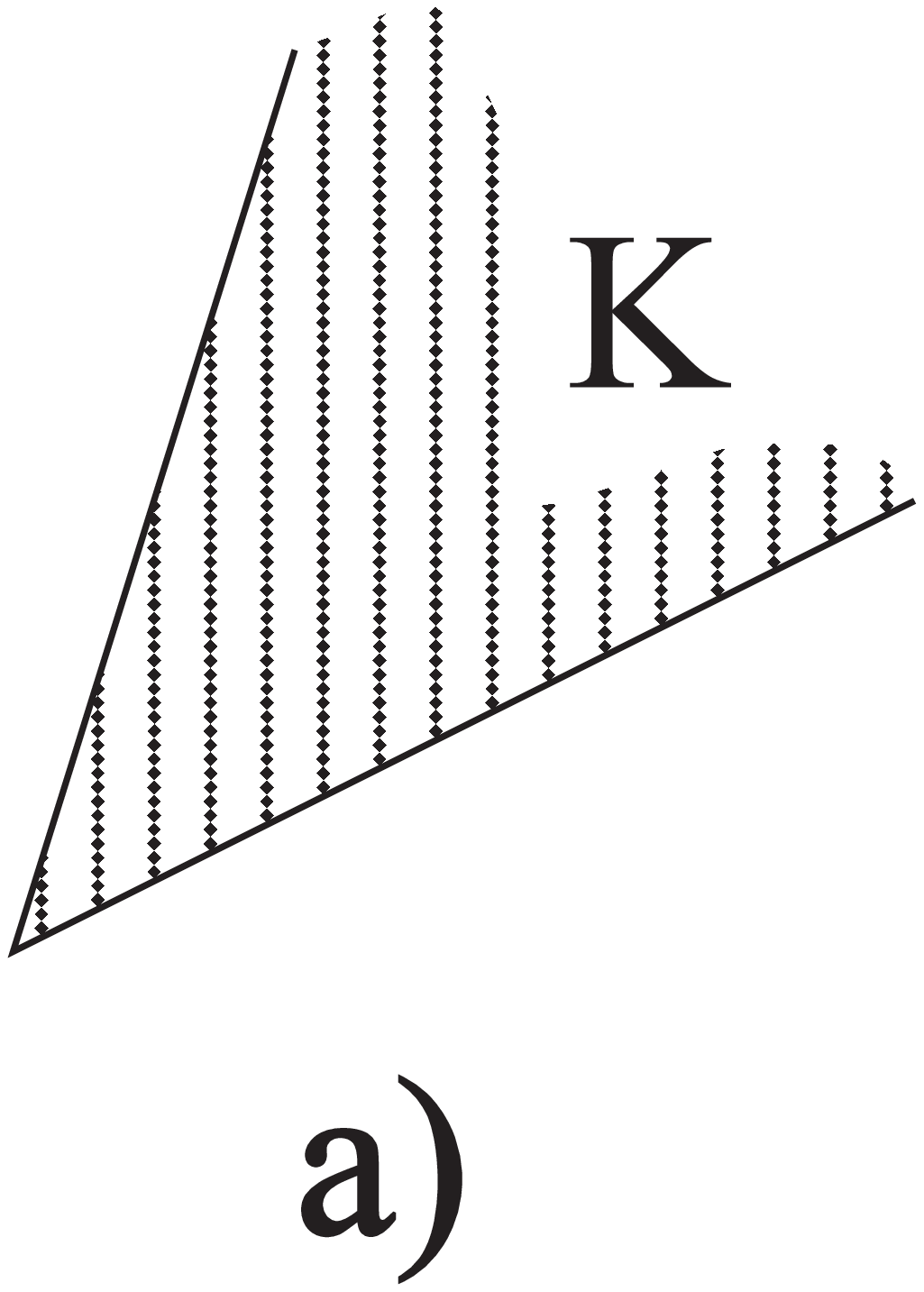}\hfill
\epsfysize=0.13\vsize\epsffile{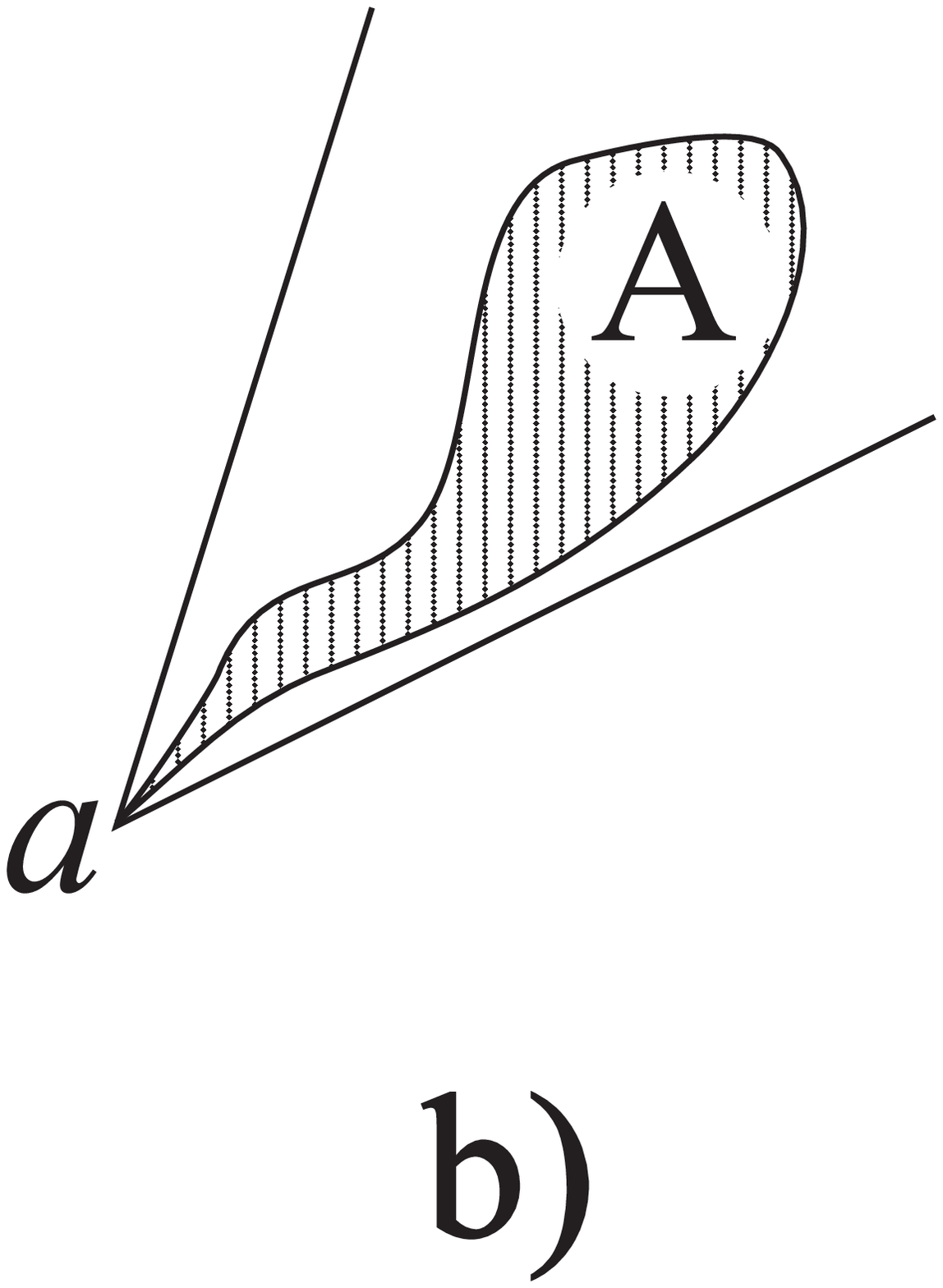}\hfill}
\botcaption{\tc Fig. 1}{\tc a) Cone $K$, b) Pointed with respect
to the cone $K$ set $A$.}
\endcaption
\endinsert

In affine space we deal with pointed cones  $K$, which are unions
of rays beginning at the vertex of the cone not containing lines.

In the projective setting it is more natural to consider cones
$\tilde K$ which are unions of lines. Such a  cone $\tilde K$ will
be called {\tci projectively pointed}, if the set of lines lying
in the cone forms a convex set in $\Bbb RP^{n-1}$. Evidently, a
cone $\tilde K$ is projectively pointed if and only if it is a
union of an affine pointed cone $K$ with its opposite cone $(-K)$,
$\tilde K=K\cup (-K)$.

We say that a set $A$ in affine space is {\tci pointed with
respect to a cone } $\tilde K=K\cup (-K)$, if the set $A$ is
pointed with respect to both the cone $K$ and the cone $(-K)$.

\midinsert
\input epsf
\centerline{\hfill\epsfysize=0.18\vsize\epsffile{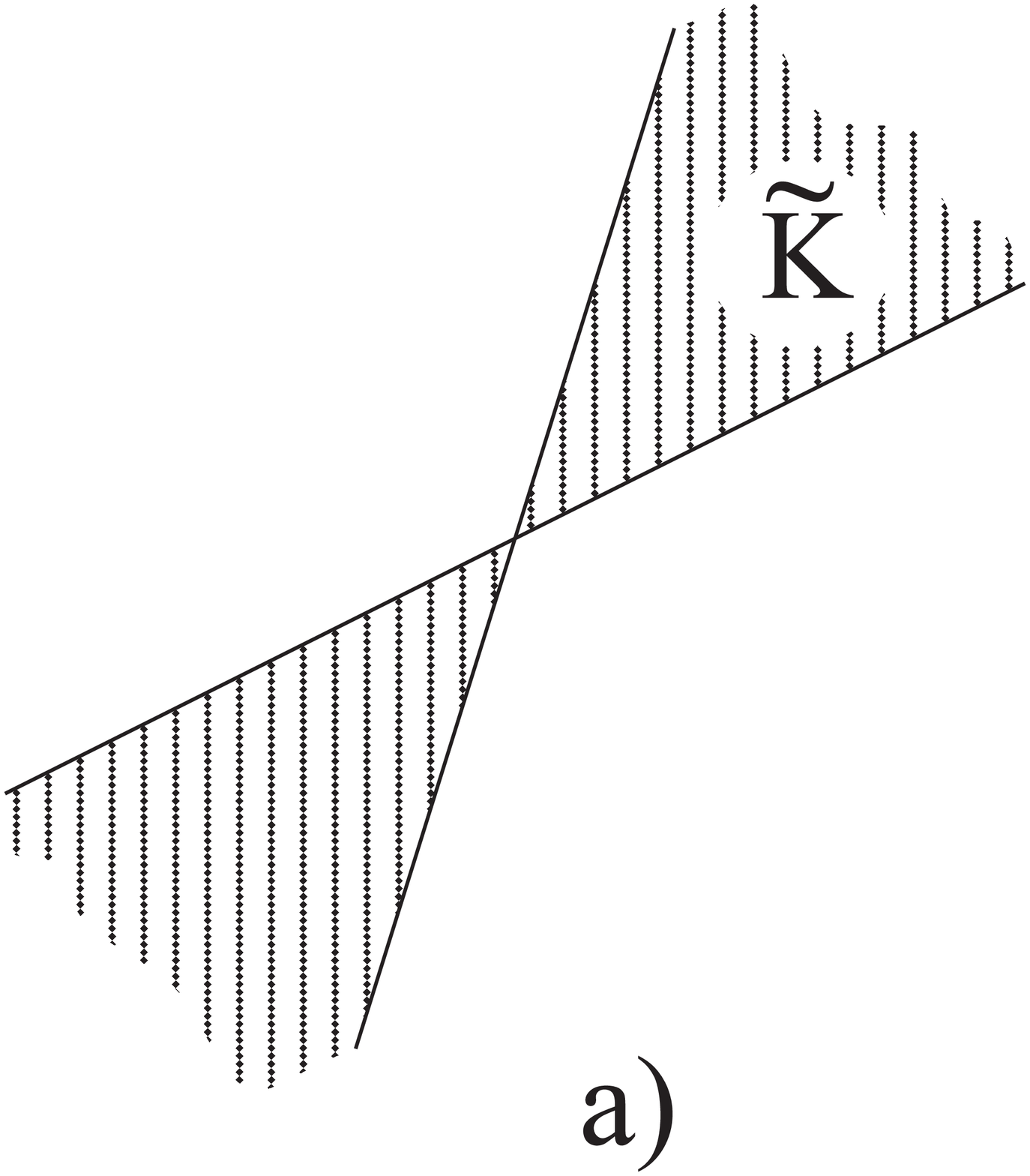}\hfill
\epsfysize=0.18\vsize\epsffile{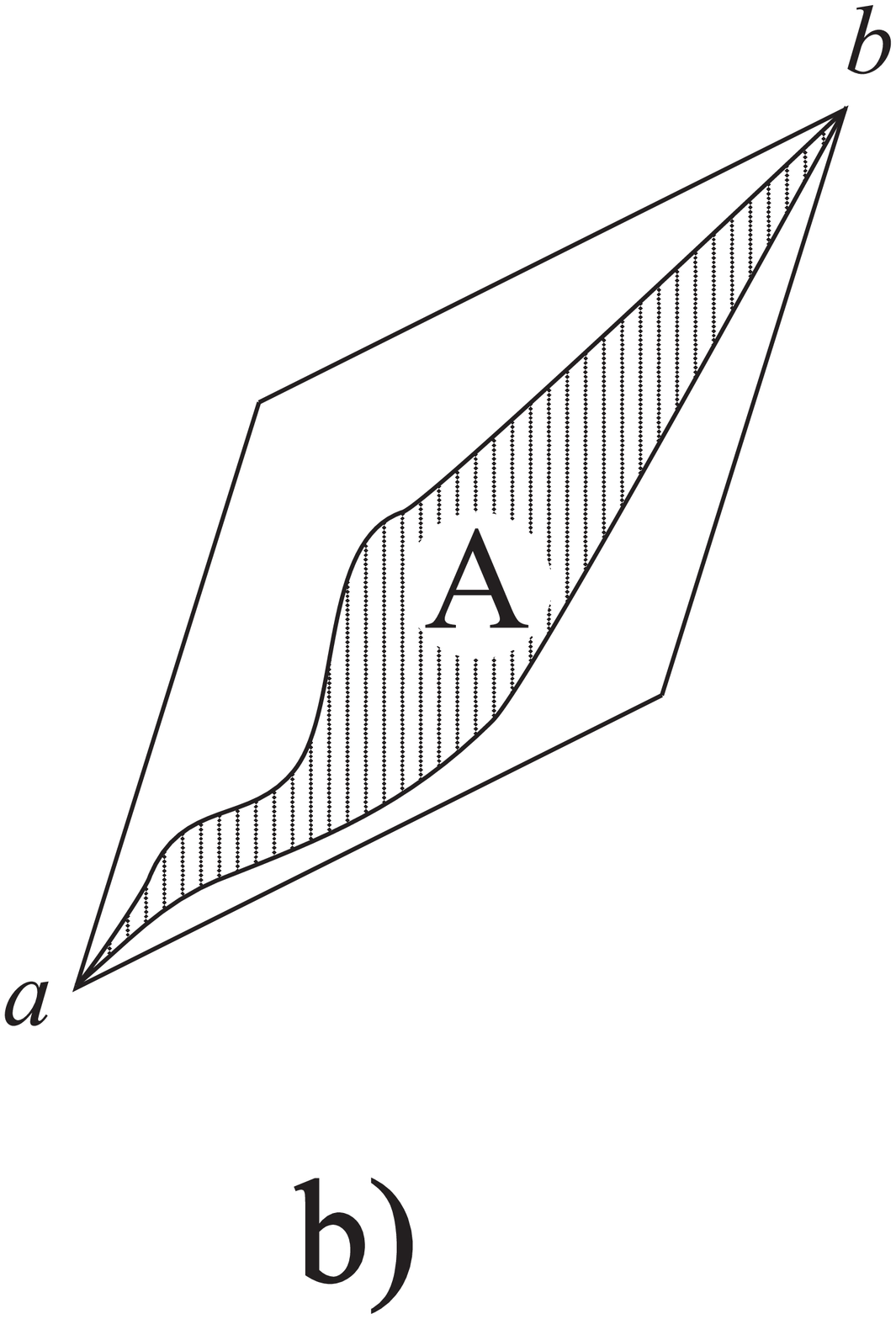}\hfill}
\botcaption{\tc Fig. 2}{\tc a) Projectively pointed cone $\tilde
K$, b)A set $A$ pointed with respect to the cone $\tilde K$.}
\endcaption
\endinsert

A set $A$ pointed with respect to a cone $\tilde K$ has two
vertices $a$ and $b$, relative to the cones $K$ and $(-K)$
correspondingly.

The following statement is evident.}

\proclaim{\tcr Proposition}{\tci Suppose  a connected set $A$ is
pointed with respect to a cone $\tilde K= K\cup (-K)$, and let $a$
and $b$ be  vertices of $A$ relative to $\tilde K$. Let $\tilde Q$
be a hyperplane intersecting $\tilde K$ at one point (the origin)
only. Then an affine hyperplane $Q$, parallel to a hyperplane
$\tilde Q$, intersects the  set $A$ if and only if $Q$ intersects
the segment joining the points $a$ and $b$. Vice versa, if a
connected set $A$ with fixed points $a$ points $b$ has this
property, then the set $A$ is pointed with respect to the cone
$\tilde K=K\cup (-K)$ and $a$ and $b$ are the vertices of $A$.}
\endproclaim

{\tc Let's turn now to a projective setting. Let $N$ be a
projective space, $L\subset N$ be a fixed hyperplane and $\Delta
\subset L$ be a closed convex set in $L$.

We say that a connected set $A\subset N$, not intersecting the
hyperplane $L$, is {\tci pointed with respect to the convex set
$\Delta$}, if there exist two points $a$ and $b$ in the set $A$
(so-called {\tci vertices} of the set $A$ with respect to
$\Delta$) such that any hyperplane $p$ in projective space $N$,
not intersecting the convex set $\Delta \subset L$, intersects $A$
if and only if $p$ intersects the segment joining the points $a$
and  $b$ and lying in the affine space $N\subset L$.

This projective definition is a projective reformulation of the
affine definition. Indeed, the projective space is a linear space
with an added hyperplane at infinity. To the convex set $\Delta$,
lying in the hyperplane at infinity, corresponds a pointed  cone
$\tilde K$ equal to the union of all lines passing through the
origin and points of the set $\Delta$.

According to the Proposition, the set $A$ in the affine space
$N\setminus L$ is pointed with respect to the cone $\tilde K$ if
and only if the set $A$, considered as a subset of  projective
space, is pointed with respect to the  convex set $\Delta=\tilde
K\cap L$.}

\subhead{\tcb Families of convex sets affinely dependent on
parameters}
\endsubhead
{\tc We begin with an affine setting. Fix a linear subspace $N$ of
a linear space $\Bbb R^n$. The linear space $\Bbb R^n$ is fibered
by affine subspaces $N_m$ parallel to $N$ and parameterized by
points $m$ of a factor-space $\Bbb R^n/N$. Fix a convex domain
$\Delta$ in the space of parameters $\Bbb R^n/N$.  Suppose that
for each point $m\in \Delta$ in the  affine space $N_m$ a closed
convex set $A_m\subset N_m$ is given.

We say that a {\tci family of convex sets $\{A_m\}$ depends
affinely on parameter} $m\in \Delta$, if for any two points $m_1,
m_2 \in \Delta$ and any $0\leq t\leq 1$, the set $A_{m_t}$
corresponding to the parameter $m_t=tm_1+(1-t)m_2$ is a linear
combination $tA_{m_1}+ (1-t)A_{m_2}$ of sets $A_{m_1}$ and
$A_{m_2}$ in Minkowski sense. }

\proclaim{\tcr Proposition}{\tci A family of convex sets $A_m$,
$m\in \Delta$ depends affinely on parameters if and only if for
any simplex  $\Delta (a_1,\dots, a_k)\subset \Delta $ with
linearly independent vertices $a_1,\dots, a_k\in \Delta$ a convex
hull of a union  of the sets $A_{a_1}\dots, A_{a_k}$ coincides
with a union of the sets $A_m$ for all parameters $m\in \Delta
(a_1,\dots,a_k)$.}
\endproclaim

{\tc A particular case of one-dimensional space $N=(l)$ is
especially simple. In this case the convex sets $A_m$ are simply
segments, and the Proposition reads as follows.}

\proclaim{\tcr Proposition}{\tci A family of parallel segments
$A_m$ in $\Bbb R^n$ depends affinely on parameter  $m$ belonging
to a convex domain $\Delta\subset \Bbb R^n/(l)$, if and only if
there exist two hyperplanes $\Gamma_1$ and $\Gamma_2$ in the space
$\Bbb R^n$ such that, first,  for any $m\in \Delta$ ends of the
segment $A(m)$ coincide with points of intersections of the line
$N_m$ with the hyperplanes $\Gamma_1$ and $\Gamma_2$, and, second,
projection along $N$ of an intersection of $\Gamma_1$ and
$\Gamma_2$ doesn't intersect the interior of $\Delta$.}
\endproclaim

{\tc The general definition of affine dependence on parameter can
be  reduced, using projections, to the  case of one-dimensional
space. Let $Q$ be a subspace of the space $N$. A quotient space
$\Bbb R^n/Q$ contains a subspace $\pi(N)=N/Q$. Spaces $(\Bbb
R^n/Q)/\pi(N)$ and $\Bbb R^n/N$ are naturally isomorphic and we
will use this isomorphism.

We say that {\tci a family of convex sets $A_m\subset N_m$ depends
affinely on parameter $m\in \Delta\subset \Bbb R^n/N$ in direction
of the hyperplane  $Q$ in space $N$}, if after the projection $\pi
\:\Bbb R^n\to \Bbb R^n/Q$ the segments $\pi(A_m)$ on lines $N_m/Q$
depend affinely on parameter $m\in \Delta$. (Using the
isomorphism of $(\Bbb R^n/Q)/\pi (N)$ and $\Bbb R^n/N$, we
consider $\Delta\subset \Bbb R^n/N$ as a set in $(\Bbb R^n/G)/\pi
(N)$.)}

\proclaim{\tcr Theorem}{\tci A family of convex sets $A_m\subset
N_m$ depends affinely on parameter $m\in \Delta$ if and only if
the family $A_m\subset N_m$ depends affinely on parameter $m\in
\Delta$ in direction of $Q$ for any hyperplane $Q$.}
\endproclaim

\demo{\tcr Proof}{\tc   Taking a subspace $M$ transversal to $N$,
we can identify all parallel spaces $N_m$ (two points of different
sections are identified if they lie in the same translate of $M$).
Then all dual spaces $N^*_m$ are identified with the space $N^*$
and all support functions $H_m(\xi)= \max\limits _{x \in A_m}(\xi,
x)$ of convex sets $A_m$ can be considered as functions on the
same space $N^*$.

To a linear combination (in Minkowski sense) of convex sets
corresponds a linear combination of their support functions. So
the dependence  of the family of convex sets $A_m$ on parameter
$m\in \Delta$ is affine if and only if for any fixed covector $\xi
\in N^*$ the support function $H_m(\xi)$ is a linear polynomial on
parameter $m$.

Let's rewrite this condition for $\xi$ and $-\xi$ simultaneously.
Denote by $Q$  a hyperplane in $N$ defined by an equation
$(\xi,x)=(-\xi,x)=0$. Project the set $A=\bigcup_{m\in \Delta}A_m$
along the space $Q$. The projection $\pi(A)$ lies in the space
$\Bbb R^n/Q$ with a marked one-dimensional subspace $l=N/Q$.
 On each line $l_m$, $m\in \Bbb R^n/N=(\Bbb R^n/Q)/l$ lies a segment
$\pi (A_m)$ equal to the projection of the convex
set $A_m$.

By assumption, the segments $\pi( A_m)$ lie between two
hyperplanes $\Gamma_1$ and $\Gamma_2$. Also, the ends $x(m)$ and
$y(m)$ of these segments lie on the line $l_m$,  and are defined
by equations $H_m(\xi)=\langle \xi,x(m)\rangle$,
$H_m(-\xi)=\langle -\xi,y(m)\rangle$. Therefore the affine
dependence of convex sets $A_m$, $m\in Q$ in direction $Q$ means
that the support functions $H_\xi (m)$ and $H_{-\xi} (m)$, where
$\xi $ are covectors orthogonal to $Q$, are polynomials of first
degree in $m\in \Delta$. Since this is true for any hyperplane
$Q\subset N$, the function  $H_\xi(m)$ depends linearly on $m$ for
any fixed $\xi$.}\enddemo

{\tc Consider now projective settings. Instead of a linear space
$\Bbb R^n$ fibered by affine subspaces $N_m$ parallel to a space
$N$ and parameterized by points of the factor-space $\Bbb R^n/N$,
we will have the following objects: a projective space $\Bbb RP^n$
with a projective subspace  $L$, fibered by subspaces $N_m$ of
dimension $\dim N_m=\dim L+1$ and containing the space $L$. The
subspaces $N_m$ are parameterized by points of a factor-space
$M=(\Bbb RP^n)/N$. Consider parameters $m$ belonging to a convex
set $\Delta \subset M$.

Let $T\subset L$ be a hyperplane in $L$. Denote a projection of
the projective space from the center $T$ by $\pi$.  Projection of
the space $L$ is just a point $\pi(L)$. Projection of the  space
$N$ is a line $l$ belonging to a bundle of all lines
$l_m=\pi(N_m)$ containing the marked point $\pi(L)$. After a
natural identification of factor-spaces $(\Bbb RP^n)/L$ and $(\Bbb
RP^n/T)/\pi(L)$, the space $N_m\subset \Bbb RP^n$ and the line
$l_m=\pi(N_m)\subset \Bbb RP^n/T$ correspond to the same parameter
$m\in \Bbb RP^n/L= (\Bbb RP^n/T)/\pi(L)$. The  domain $\Delta
\subset \Bbb RP^n/L$ can be considered as a domain in the space
$(\Bbb RP^n/T)/\pi(L)$.

Introduce the following notation. Let $\Gamma _1$ and
$\Gamma _2$ be two hyperplanes in projective space, not containing
the point $\pi(L)$, and $l$ be a line containing this point.
Points of intersection of $\Gamma _1$ and $\Gamma_2$ with the line
$l$ divide it into two segments. The segment not containing the
point $\pi(L)$ will be called exterior relative to the point $\pi
(L)$ segment between hyperplanes $\Gamma _1$ and $\Gamma _2$ on
the line $l$.

Let $A$ be a set  not intersecting space $L$, whose sections $A_m$
by the spaces $N_m\supset L$ are convex. We say that {\tci
sections $A_m$ depend affinely on parameter}  $m$ belonging to a
convex domain $\Delta \subset \Bbb RP^n/L$ {\tci in direction of
the hyperplane} $T\subset L$, if the sections of the set $\pi(A)$
by lines $l_m$ containing the point $\pi(L)$, depend affinely on
$m\in\Delta\subset \Bbb RP^n/L(=(\Bbb RP^n/T)/(\pi(L))$. In other
words, there exist two hyperplanes $\Gamma _1$ and $\Gamma _2$ in
$\Bbb RP^n/T$, not containing $\pi(L)$, such that, first, the
intersection of $\pi (A)$ with any line $l_m$, $m\in \Delta$, is
equal to the exterior relative to $\pi(L)$ segment of the line
$l_m$ lying between $\Gamma _1$ and $\Gamma _2$, and, second, the
projection of $\Gamma_1\cap\Gamma_2$ on $\Bbb RP^n/L$ doesn't
intersect $\Delta$.

Now we can give a definition of affine dependence of sections on
parameter belonging to a convex  domain of the space of
parameters.

We say that {\tci sections $A_m$ of the set} $A$ by projective
spaces $N_m\supset L$ {\tci depend affinely on parameter} $m$ in
domain $\Delta$, if $A_m$ depend affinely on parameter $m$ in
domain $\Delta$ with respect to any hyperplane $T\subset L$. The
following statement can be easily checked.}

\proclaim{\tcr Proposition}{\tci Let $\Gamma $ be a projective
hyperplane containing the space $L$, such that its projection to
the space $(\Bbb RP^n)/L$ doesn't intersect a convex set $\Delta
\subset \Bbb RP^n/L$. Consider an affine chart $U$ of the
projective space, $U=\Bbb RP^n\setminus \Gamma $. Sections of a
set $A\subset \Bbb RP^n$, $A\cap L=\emptyset$, by spaces $N_m$
depend affinely on  parameter $m$ in domain $\Delta$, if and only
if the sections of the set $A\cap U$ in the affine space $U$ by
parallel spaces $N_m\setminus \Gamma $ depend affinely on
parameter $m$ in domain $\Delta\subset ( (\Bbb RP^n)/L )\setminus
(\Gamma /L)$.}
\endproclaim

\subhead{\tcb Duality}
\endsubhead
{\tc Let $\Delta $ be a convex domain in the space $L$, and let
$\Delta ^*_p$ be a dual convex domain in the space $(\Bbb
RP^n)^*/L^*$. The space $(\Bbb RP^n)^*/L^*$ parameterizes $(\dim
L^*+1)$-dimensional subspaces of $(\Bbb RP^n)^*$ containing $L^*$.
The domain $\Delta^*$ corresponds  to subspaces $Q^*\subset(\Bbb
RP^n)^*$ of this type dual to subspaces $Q\subset L$ not
intersecting the domain $\Delta$.}

\proclaim{\tcr Theorem}{\tci Let $A$ be a $L$-convex-concave
subset of a projective space $\Bbb RP^n$. A section $A\cap N$ of
the set $A$ by a $(\dim L+1)$-dimensional subspace $N$ containing
$L$, is pointed relative to a convex domain $\Delta\subset L$, if
and only if the following dual condition holds: $L$-dual to the
$A$ subset $A^\bot_L$ of the dual space $(\Bbb RP^n)^*$ depends
affinely on parameters belonging to the domain $\Delta ^*_p$ in
direction of the hyperplane $N^*\subset L^*$.}
\endproclaim

\demo {\tcr Proof}{\tc The set $A$ is $L$-convex-concave, so the
section $A\cap N$ is dual to the complement to the projection from
the center  $N^*\subset L^*$ of the set $A^\bot _L$.

Let $a$ and $b$ be vertices of the pointed set $A\cap N$ relative
to the convex set $\Delta\subset L$. Fix a
 hyperplane $q_L$ in the space $L$, not intersecting the convex set
 $\Delta \subset L$. Consider a one-dimensional bundle  $\{p^t\}$ of hyperplanes
containing the space $q_L$ in space $N$. This bundle contains the
following three hyperplanes: the hyperplane $L$, a  hyperplane
$p_a$ containing the vertex $a$ of the set $A$, and hyperplane
$p_b$, containing the vertex $b$ of the set $A$.

Take a segment $[p_a,p_b]$ with ends corresponding to $p_a$ and
$p_b$ and not containing the point $L$ on a projective line
corresponding to the one-dimensional bundle of hyperplanes
$\{p^t\}$. Any hyperplane $p^t$ (except the hyperplane $L$ itself)
intersects $L$ by a subspace $q_L$, and $q_L$ doesn't intersect
$\Delta$. The set $A$ is pointed with respect to $\Delta$, so a
hyperplane $p^{t_0}$ intersects $A\cap N$ if and only if the point
$p^{t_0}$ belongs to the segment $[p_a,p_b]$.

Consider the dual space $(\Bbb RP^n)^*$. To the section $A\cap N$
of the set $A$ corresponds a projection of the set $A^\bot _L$
from the center $N^*$. To hyperplanes in $N$ correspond points in
a factor-space $(\Bbb RP^n)^*/N^*$. In particular, to the
hyperplane $L$ in $N$ corresponds the marked point $\pi(L^*)$ in
the factor-space $(\Bbb RP^n)^*/N^*$, namely the projection of the
space $L^*$ from the center $N^*$. To the bundle of hyperplanes
$\{p^t\}$ corresponds a line passing through $\pi(L^*)$. This line
intersects projection of the set $A^\bot _L$ exactly by a segment
$[p^a,p^b]$ not containing the point $\pi(L^*)$.

To different hyperplanes $q_L$ in space $L$ correspond different
one-dimensional bundles of hyperplanes $\{p^t\}$ in $N$, i.e.
different lines in $(\Bbb RP^n)^*/N^*$, containing the marked
point $\pi(L^*)$. The hyperplane $q_L$ in space $L$ doesn't
intersect set $\Delta$, so the dual space $q^*_L\supset L^*$ is
parameterized by a point of $\Delta_p^*\subset (\Bbb RP^n)^*/L^*$.
Projection of the space $q^*_L$ from the center $N^*$ is a line in
the space $(\Bbb RP^n)^*/N$, parameterized by the same point of
the domain $\Delta^*$. Each such line intersects projection of the
set $A^\bot_L$ by segment $[p^a,p^b]$. The point $p^a$ lies in the
hyperplane $\Gamma ^a$ in the space $(\Bbb RP^n)^*/N^*$ dual to
the point $a\in N$. The point $p^b$ lies in the hyperplane $\Gamma
^b$ in the space $(\Bbb RP^n)^*/N^*$ dual to the point $b\in N$.

Two hyperplanes $\Gamma ^a$ and $\Gamma ^b$ divide the  space
$(\Bbb RP^n)^*/N^*$ into two parts. Denote by $\Gamma (a,b)$ the
closure of the part not containing the point $\pi(L^*)$. We just
proved that the set $\Gamma (a,b)$ and projection of the set
$A^\bot _L$ to the space $(\Bbb RP^n)^*/N^*$ have the same
intersections with lines passing through the point $\pi(L^*)$, and
parameterized by point of the domain $\Delta^*_p$. The theorem
proved.}
\enddemo

\head{\tcb \S5.  $L$-convex-concave sets with planar sections
being octagons with four pairs of parallel sides}
\endhead
{\tc Consider a subset $A$ of $\Bbb RP^n$ convex-concave with
respect to a one-dimensional space $L$.  Fix four points
$a_1,\dots, a_4$ lying on the line  $L$ in this order. These
points divide $L$ into four pairwise non-intersecting intervals
$\langle a_1 , a_2 \rangle$, $\langle a_2 , a_3 \rangle$, $\langle
a_3 , a_4 \rangle$, $\langle a_4 , a_1 \rangle$. Denote their
complements to $L$ by $I_1=[a_1,a_2]=L\setminus \langle
a_1,a_2\rangle,\dots, I_4=[a_4,a_1]=L\setminus \langle
a_4,a_1\rangle$  (these segments are intersecting). In this
paragraph we prove the Main  Hypothesis for $L$-convex-concave
sets $A$ whose sections $N\cap A$  by two-dimensional planes  $N$
containing the line  $L$, are pointed relative to the segments
$I_1,\dots,I_4$.}

\proclaim{\tcr Theorem}{\tci Suppose that all planar sections
$A\cap N$ of a $L$-convex-concave set $A$, $A\subset \Bbb RP^n$,
$\dim L=1$, by two-dimensional planes $N$ containing $L$, are
pointed with respect to four segments $I_1,\dots, I_4$ on the line
$L$. Suppose that the union of $I_i$ coincides with $L$  and that
the complements $L\setminus I_i$ are pairwise non-intersecting.
Then the set $A$ contains a projective space $M$ of dimension
$(n-2)$.}
\endproclaim

{\tc Before the proof we will make two remarks.

First, the assumptions of the theorem about the convex-concave set
$A$, are easier to understand in an affine chart $\Bbb R^n$ not
containing the line $L$. In this chart the family of
two-dimensional planes containing $L$ becomes a family of parallel
two-dimensional planes. In the space $\Bbb R^n$ four classes of
parallel lines are fixed, each passing through one of the points
$a_1,\dots, a_4$ of the line $L$ at infinity. The assumptions of
the theorem mean that each section of the set $A$ by a plane $N$
is an octagon with sides belonging to these four fixed classes of
parallel lines. (Some sides of this octagon can degenerate to a
point, and number of sides of the octagon $(A\cap N)$ will then be
smaller than 8.)

Also, there is a natural isomorphism between $(\Bbb RP^1)^*$ and
$\Bbb RP^1$. Indeed, each point $c\in \Bbb RP^1$ of a projective
line is also a hyperplane in $\Bbb RP^1$. However, a segment
$[a,b]$ on the projective line $\Bbb RP^1$ will be dual to its
{\it complement} $\langle a,b\rangle= (\Bbb RP^1)\setminus [a,b]$,
and not to itself. Indeed, by definition, a dual to a convex set
$\Delta$ set $\Delta ^*_p$, consists of all hyperplanes not
intersecting $\Delta$. }

\demo{\tcr Proof of the Theorem}{\tc Consider the dual projective
space $(\Bbb RP^n)^*$ and its subspace $L^*$, $\dim L^*=(n-2)$,
dual to the line $L$. Projective line  $(\Bbb RP^n)^*/L^*$,
isomorphic to a line dual to $L$, is divided by points
$a_1^*,\dots, a_4^*$ into four intervals $\langle a_1^*,
a_2^*\rangle$, $\langle a_2^*, a_3^*\rangle$, $\langle a_3^*,
a_4^*\rangle$, $\langle a_4^*, a_1^*\rangle$, dual to segments
$I_1,\dots,I_4$. The set $A^\bot_L$, $L$-dual to the set $A$, will
depend affinely on parameter on these intervals, since the set $A$
is pointed relative to the segments $I_1,\dots ,I_4$. Therefore
the set $A^\bot _L$ is a linear interpolation of its four
sections. In other words, this set has four sections by planes
corresponding to $a_1^*, \dots, a_4^*$, and all other sections of
$A^\bot _L$ are affine combinations (in Minkowski sense) of
sections corresponding to the ends of intervals.
$L^*$-convex-concave sets of this type contain a line (see
Introduction and our paper {\it ``A convex-concave domain in
${\Bbb RP}^3$ contains a line''}, in preparation). Denote this
line by $l$. The set $A$ will contain an $(n-2)$-dimensional space
$l^*\subset \Bbb RP^n$ dual to the line $l$, q.e.d.}
\enddemo

\head{\tcb \S6. Surgeries on  convex-concave sets}
\endhead

{\tc In this section we describe two special surgeries on
$L$-convex-concave subsets of $\Bbb RP^n$, one applicable when
$\dim L=n-2$ and another when $\dim L=1$. These two surgeries are
dual. }

\subhead{\tcb The first surgery: $\dim L=n-2$}
\endsubhead
{\tc To a $(n-2)$-dimensional subspace $L$ of $\Bbb RP^n$
corresponds  a one-dimensional bundle of hyperplanes containing
$L$. These hyperplanes are parameterized by points of the
projective line $\Bbb RP^n/L$. Fix two points $a$ and $b$ and a
segment $[a,b]$ on this line , one of two segments into which the
points $a$ and $b$ divide the projective line $\Bbb RP^n/L$.

For  any $L$-convex-concave set $A$ and the segment $[a,b]\subset
\Bbb RP^n/L$ we define a set $S_{[a,b]}(A)$, which is also
$L$-convex-concave. Here is the definition of the set
$S_{[a,b]}(A)$. The hyperplanes $\Gamma_a$ and $\Gamma_b$
corresponding to parameters $a$ and $b$, $L= \Gamma_a\cap
\Gamma_b$, divide the set $\Bbb RP^n\setminus L$ into two
half-spaces: the first half-space $\Gamma ^1[a,b]$ is projected to
the segment $[a,b]$, and the second one $\Gamma ^2[a,b]$ is
projected to its complement.

Let $c$ be some point on the line $\Bbb RP^n/L$ not belonging to
the segment $[a,b]$, and let $\Gamma _c$ be the corresponding
hyperplane in $\Bbb RP^n$.}

\definition{\tcr Definition}{\tc The set $S_{[a,b]}(A)$
is defined by the following requirements:
\item{1)} the set $S_{[a,b]}(A)$ doesn't intersect the space $L$,
i.e. $S_{[a,b]}(A)\cap L= \emptyset$;
\item{2)} the set $S_{[a,b]}(A)\cap
\Gamma ^1_{[a,b]}$ coincides with  a convex hull of the union of
sections $A\cap \Gamma (a) $ and $A\cap \Gamma (b)$ in an affine
chart $\Bbb RP^n\setminus \Gamma _c$;
 \item{3)} the set  $S_{[a,b]}(A)\cap \Gamma ^2_{[a,b]}$ coincides with  $A\cap \Gamma ^2_{[a,b]}$.
}
\enddefinition

{\tc It is easy to see that the set $S_{[a,b]}(A)$ is correctly
defined,  i.e. it doesn't depend on the choice of the hyperplane $\Gamma _c$.}

\proclaim{\tcr Theorem}{\tci For any $L$-convex-concave set $A$
the set $S_{[a,b]}(A)$ is also $L$-convex-concave.}
\endproclaim

\demo{\tcr Proof}{\tc Any section of the set $S_{[a,b]}(A)$ by a
hyperplane $\Gamma _d$ containing $L$, is convex. Indeed, if
$d\notin [a,b]$, then $\Gamma _d \cap S_{[a,b]}(A)=\Gamma _d\cap
A$, and the  set $\Gamma _d\cap A$ is convex by definition.
Otherwise, i.e. if $d\in [a,b]$, the  $\Gamma _d\cap S_{[a,b]}(A)$
is a linear combination (in Minkowski sense) of convex sections
$A\cap \Gamma _a$ and $A\cap \Gamma _b$ (in any affine chart $\Bbb
RP^n\setminus \Gamma _c$, $c\notin [a,d]$), so is convex .

Let's prove that  a complement to a projection of the set
$S_{[a,b]}(A)$ from any center $L_1\subset L$, where $L_1$ is a
hyperplane in $L$, is a convex open set. Consider a projection
$\pi(A)$ of the set $A$ on the projective plane $\Bbb RP^n/L_1$.
The set $A$ is $L$-convex-concave, so the complement $B$ to the
projection $\pi (A)$ is an open convex set, containing the marked
point $\pi (L)$. The plane $\Bbb RP^n/L_1$ contains  two lines,
$l_a=\pi (\Gamma _a)$ and $l_b=\pi (\Gamma _b)$, passing through
the point $\pi (L)$, a half-plane $l^1_{[a,b]}=\pi (\Gamma
^1_{[a,b]})$ and a  complementary half-plane
$l^2_{[a,b]}=\pi(\Gamma ^2_{[a,b]})$.

From the definition of the set $S_{[a,b]}(A)$ we see that the
complement $B_{[a,b]}$   to its projection $\pi (S_{[a,b]}(A))$
has the following structure.
\item{1)} set  $B_{[a,b]}$ contains the point $\pi(L)$;
\item{2)}  Consider two closed
triangles with vertices at the point  $\pi(L)$, lying in
$l^1[a,b]$,  with one side being the segment lying inside
$l^1[a,b]$ and joining the points of intersection of lines $l_a$
and $l_b$ with the  boundary of the domain $B$ (see Fig.~3). The
set  $B_{[a,b]}\cap l^1[a,b]$ is a union of these triangles with
the sides described above removed;
\item{3)} set  $B_{[a,b]}\cap l^2[a,b]$ coincides with the set $B\cap
l^2[a,b]$.

From this description we see that the set $B_{[a,b]}$ is convex
and open, q.e.d.}
\enddemo

\midinsert
\input epsf
\centerline{\hfill\epsfysize=0.18\vsize\epsffile{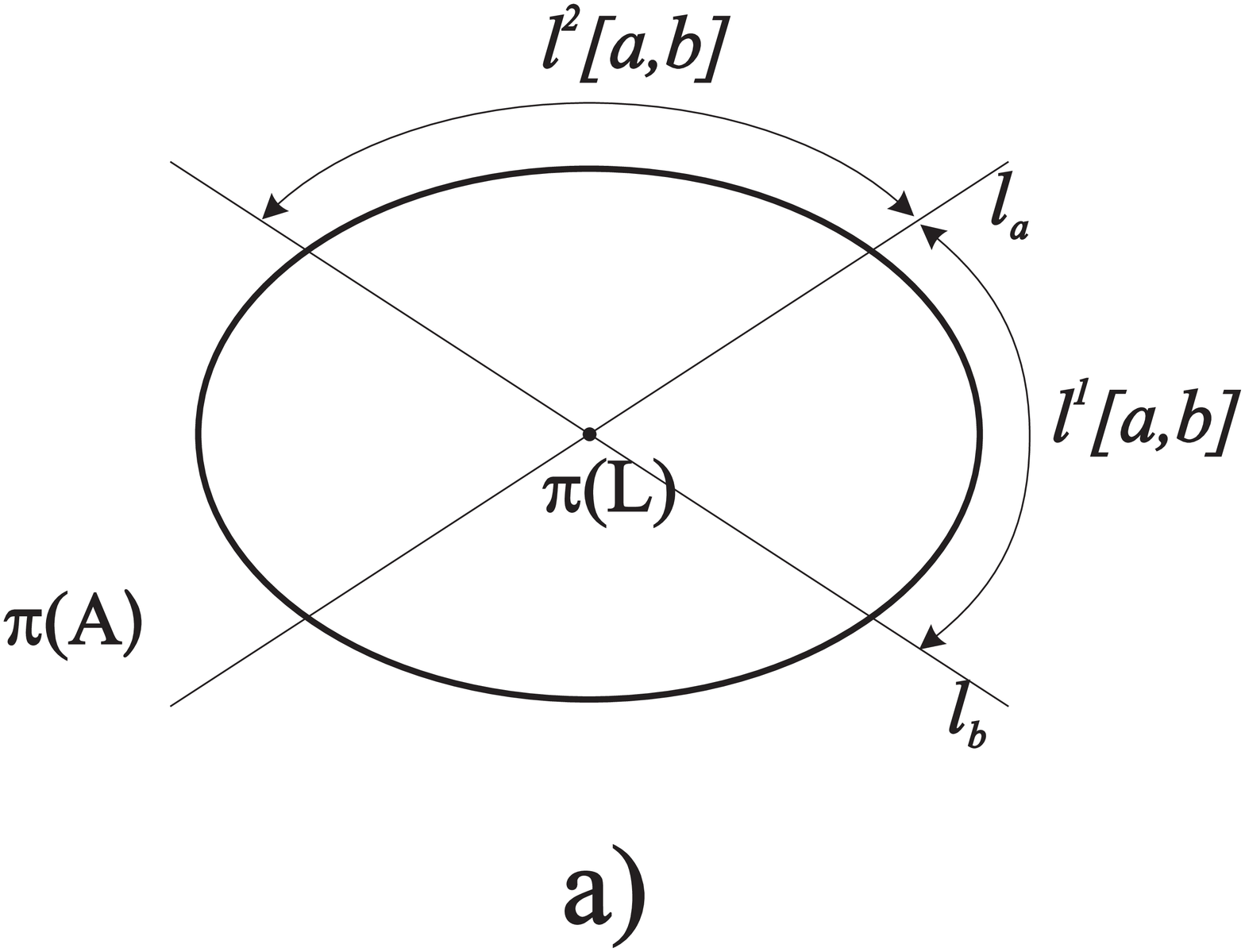}\hfill
\epsfysize=0.18\vsize\epsffile{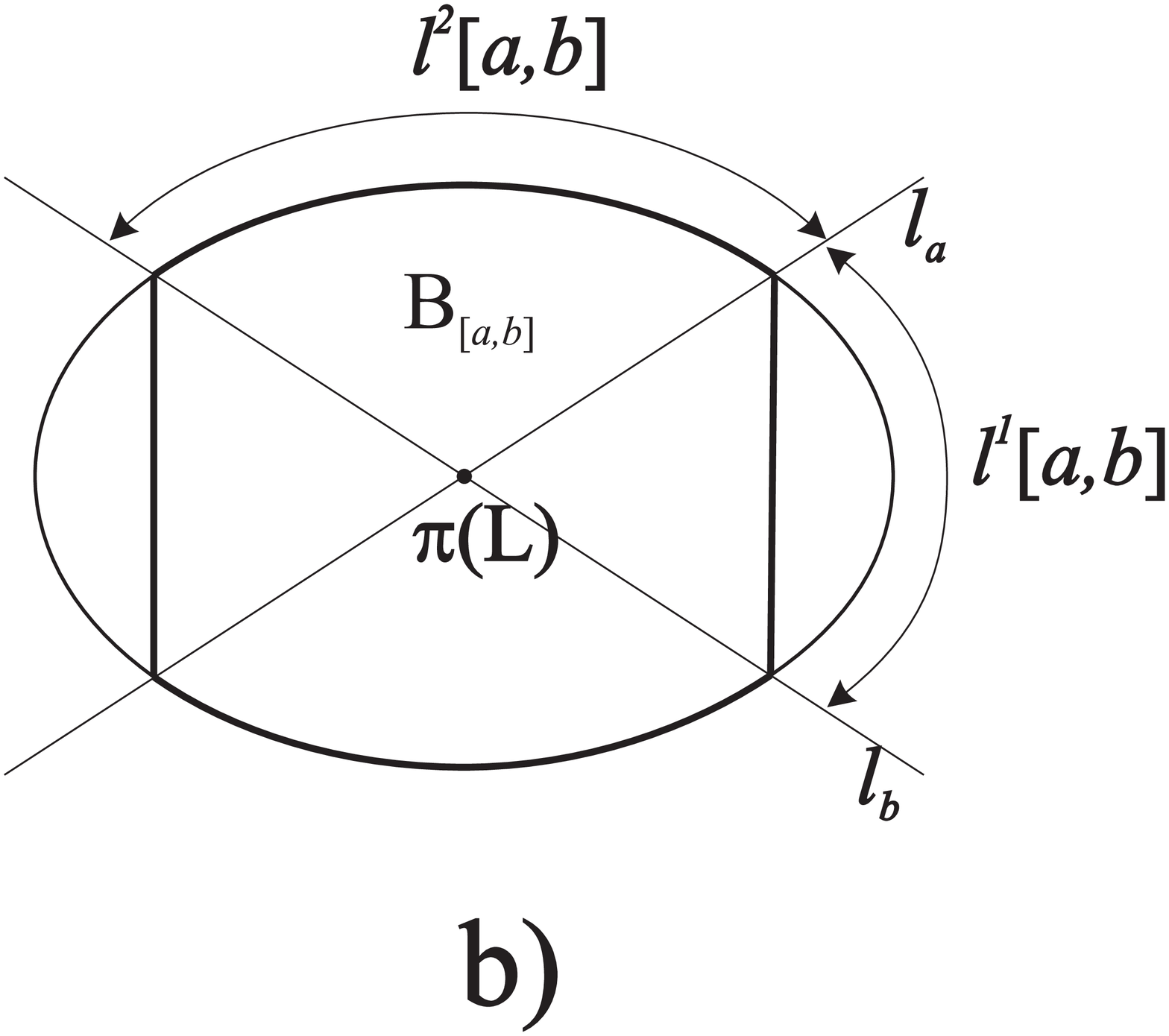}\hfill}
\botcaption{\tc Fig. 3}{\tc a) The complement  $B$ to the
projection $\pi (A)$ of the set $A$, b) The complement $B_{[a,b]}$
to the projection $\pi (S_{[a,b]}(A))$ of the set $S_{[a,b]}(A)$.}
\endcaption \endinsert

{\tc If two segment $[a,b]$ and $[c,d]$ on the line $\Bbb RP^n/L$
do not have  common interior points, then the surgeries
$S_{[a,b]}$ and $S_{[c,d]}$ commute. We can   divide the line
$\Bbb RP^n/L$  into a finite set of segments $[a_1,a_2],\dots,
[a_{k-1}, a_{k}], [a_{k},a_1]$ and apply to a $L$-convex-concave
set $A$ the surgeries corresponding to these segments. As a result
we will get a $L$-convex-concave set $D$, such that the  sections
of $D$ by hyperplanes $\Gamma _{a_1},\dots, \Gamma _{a_n}$
coincide with  sections $A\cap \Gamma_{a_i}$ of the set $A$ by the
same hyperplanes. For an intermediate value of parameter
$a_i<a<a_{i+1}$ the section $D\cap \Gamma_a$ coincides with the
section by the same hyperplane of the convex hull of the union of
sections $A\cap \Gamma _{a_i}$ and $A\cap \Gamma _{a_{l+1}}$ in
affine chart $\Bbb RP^n\setminus \Gamma _c$ (where $c$ is any
point of the line $\Bbb RP^n/L$, not belonging to the segment
$[a_i,a_{i+1}]$).}

\subhead{\tcb The second surgery: $\dim L=1$}
\endsubhead
{\tc To a  one-dimensional space  $L$ corresponds a
$(n-2)$-dimensional bundle of two-dimensional planes containing
the line $L$. Fix two points $a$ and $b$ and a segment $[a,b]$ on
the line $L$ --- one of two segments into which the points $a$ and
$b$ divide the line $L$. For any $L$-convex-concave set $A$ and
the segment $[a,b]\subset L$ we construct a new $L$-convex-concave
set $P_{[a,b]}(A)$. A section of $P_{[a,b]}(A)$ by any
two-dimensional plane $N$, $N\supset L$, will depend on the
section of the set $A$ by this plane $N$ only.

We define first an operation $F_{[a,b]}$ applicable to
two-dimensional convex sets. This operation $F_{[a,b]}$ transforms
planar sections $A\cap N$ of the set $A$ to  planar sections
$P_{[a,b]}(A)\cap N$ of the set $P_{[a,b]}(A)$.

Consider a two-dimensional projective plane $N$ with a
distinguished projective line $L$ and a segment $[a,b]\subset L$.
Let $\Delta \subset N$ be any closed convex subset of the plane
$N$, not intersecting the line $L$.

By definition, the operation  $F_{[a,b]}$ transforms a set $\Delta
\subset N$ to the smallest convex set $F_{[a,b]}(\Delta)$
containing the set $\Delta$ and pointed relative to the segment
$[a,b]$.

Here is a more explicit description of the set
$F_{(a,b)}(\Delta)$.

\midinsert
\input epsf
\centerline{\hfill\epsfysize=0.18\vsize\epsffile{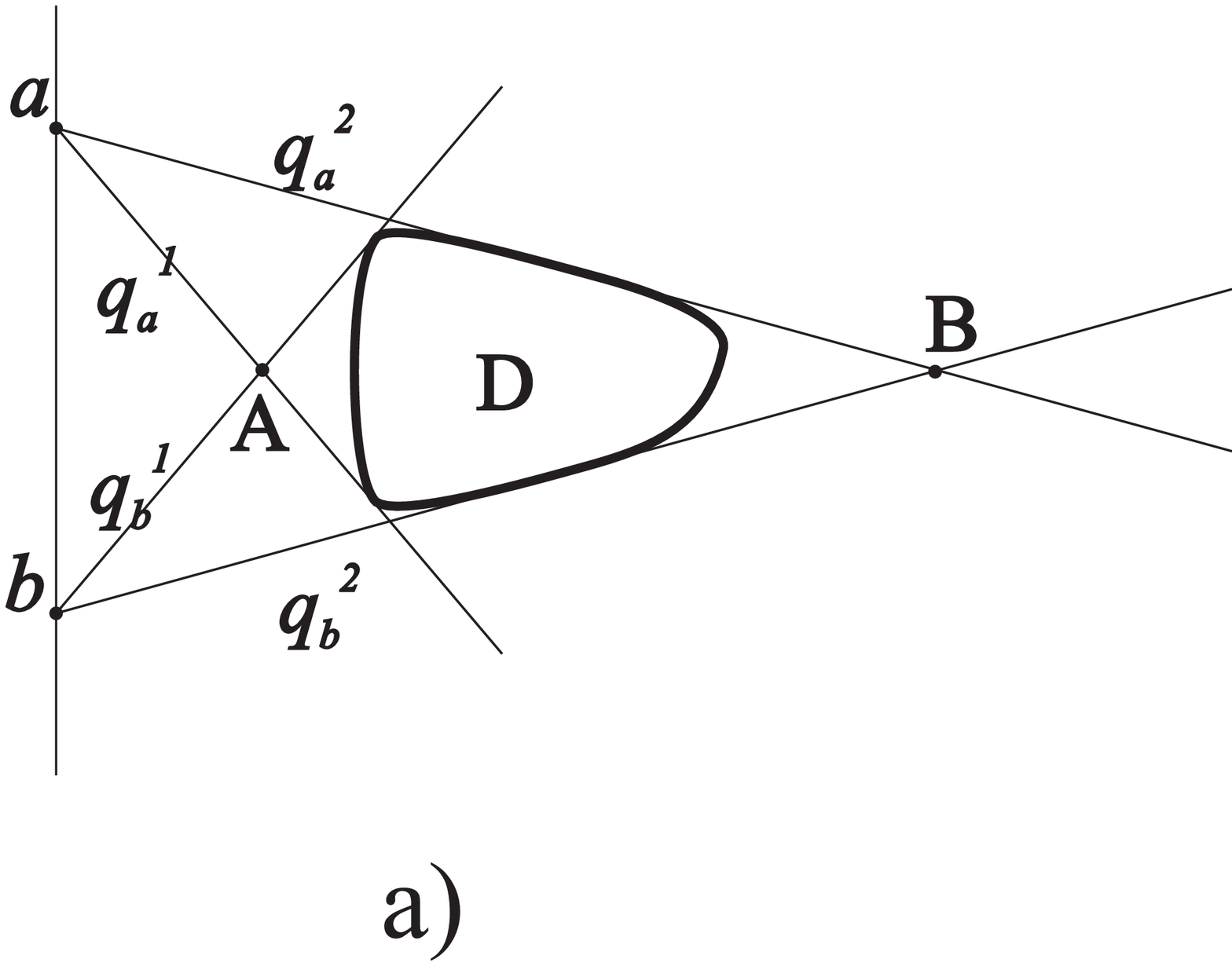}\hfill
\epsfysize=0.18\vsize\epsffile{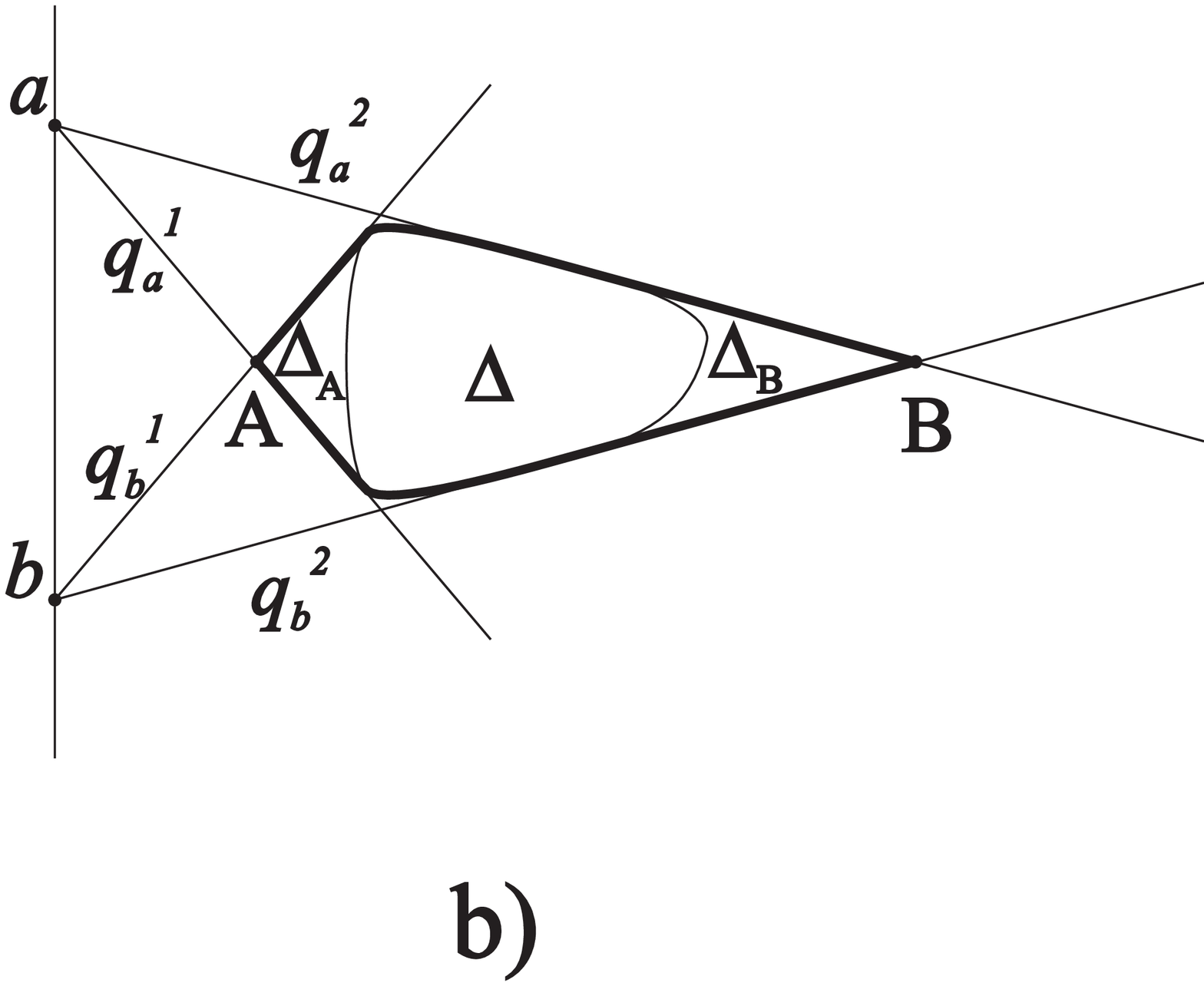}\hfill}
\botcaption{\tc Fig. 4}{\tc a) The set $\Delta$ and tangents to it
passing through the points $a$ and $b$, b) The set
$F_{(a,b)}(\Delta)$.}
\endcaption
\endinsert

Draw four tangents,  $q^1_a$, $q^2_a$ and $q^1_b$, $q^2_b$, to the
set  $\Delta$ passing through the points $a$ and $b$
correspondingly (see Fig.~4). In the convex quadrangle $\Delta _1$
in the affine plane $N\setminus L$ with sides on the lines
$q^1_a$, $q^2_a$ and $q^1_b$, $q^2_b$ there are exactly two
vertices $A$ and $B$ satisfying the following condition: the
support lines to the quadrangle $\Delta _1$ at the vertex do not
intersect the segment $[a,b]$. To the vertex $A$ corresponds a
curvilinear triangle $\Delta _A$ with two sides lying on two sides
of the quadrangle $\Delta _1$ joint to the vertex $A$. The third
side of $\Delta _A$  coincides with  the part of the boundary of
the set $\Delta $ visible from the point  $A$.

A similar curvilinear triangle $\Delta _b$ corresponds to the
vertex $B$. Evidently the set $F_{[a,b]}(\Delta)$ coincides with
the set $\Delta _A\cup \Delta \cup \Delta_B$.

Now we can define the set $P_{[a,b]}(A)$.}

\definition{\tcr Definition}{\tc For any $L$-convex-concave
subset $A$ of $\Bbb RP^n$, $\dim L=1$, and for any segment $[a,b]$
of the line $L$ we define the set $P_{[a,b]}(A)$ by the following
condition: a section $P_{[a,b]}(A)\cap N$ of this set by any
two-dimensional plane $N$ containing $L$  is obtained from the
 section $A\cap N$ of the set $A$ by operation
$F_{[a,b]}$ in the plane $N$: $P_{[a,b]}(A)\cap N=F_{[a,b]}(A\cap
N)$.}
\enddefinition

\proclaim{\tcr Theorem}{\tci For any $L$-convex-concave set $A$,
$\dim L=1$, and any segment $[a,b]\subset L$ on the line $L$ the
set $P_{[a,b]}(A)$ is $L$-convex-concave.}
\endproclaim

\demo{\tcr Proof}{\tc To any  $L$-convex-concave set $A$ in $\Bbb
RP^n$ corresponds its $L$-dual $D=(A^\bot_L)$ in the dual
projective space $(\Bbb RP^n)^*$. The set $D$ is a
$L^*$-convex-concave set, and $\dim L^*=n-2$.
 The line $L$ is dual to the set of parameters $(\Bbb
RP^n)^*/L^*$. To the segment $[a,b]\subset L$ corresponds a dual
interval $\langle a^*,b^*\rangle \subset (\Bbb RP^n)^*/L^*$. By
the segment $[a^*,b^*]\subset (\Bbb RP^n)^*/L^*$ and the
$L^*$-convex-concave set $D=A^\bot _L$ we define a new
$L^*$-convex-concave set $S_{[a^*,b^*]}(D)$. To prove the theorem
it is enough to check that the set $P_{[a,b]}(A)$ is $L^*$-dual to
the set  $S_{[a^*,b^*]}(D)$, where $D=A^\bot _L$. This is
proved below. }\enddemo

\proclaim{\tcr Proposition}{\tci The set $P_{[a,b]}(A)$ is
$L^*$-dual to the  set $S_{[a^*,b^*]}(D)$.}
\endproclaim

\demo{\tcr Proof}{\tc We proved above that if the set $D$ is
$L^*$-convex-concave, then the set $S_{[a^*,b^*]}(D)$ is also
$L^*$-convex-concave and described how to obtain the planar
projections of the set  $S_{[a^*,b^*]}(D)$ from the planar
projections of the set $D$.

Consider the sets $D^\bot _{L^*}=A$ and $S_{[a^*,b^*]}(D)^\bot
_{L^*}$ $L^*$-dual to the sets  $D$ and $S_{[a^*,b^*]}(D)$
correspondingly. Planar projections of the sets $D$ and
$S_{[a^*,b^*]}(D)$ are dual to the planar sections of the sets $A$
and $(S_{[a^*,b^*]}(D))^\bot_{L^*}$. Looking on the planar
pictures, one easily sees that sections of the set
$(S_{[a^*,b^*]}(D))^\bot_{L^*}$ are obtained from the sections of
the set $A$ by the surgery $F_{[a,b]}$. Therefore
$(S_{[a^*,b^*]}(D))^\bot _{L^*}= P_{[a,b]}(A)$.}
\enddemo

{\tc If the complements $\langle a,b\rangle_0$ and $\langle
c,d\rangle _0$ to the segments $[a,b]$ and $[c,d]$ do not
intersect, then the operations $P_{[a,b]}$ and $P_{[c,d]}$
commute. Divide the line $L$ into a finite number of intervals
$\langle a_1,a_2\rangle_0, \dots , \langle a_{k+1},a_1\rangle_0$,
complementary to segments $[a_1, a_2],\dots,[a_{k+1}, a_1]$ (the
segments intersect one another) and apply to the
$L$-convex-concave set $A$ the operations $P_{[a_i,a_{i+1}]}(A)$
corresponding to all these segments. As a result we will get a
$L$-convex-concave set $D$, whose  section by any two-dimensional
plane $N$ containing the line $L$, is a polygon with $2k$ sides
circumscribed around the section $A\cap N$ (some of the sides of
the resulting polygons can degenerate into points). To each point
$a_i$ correspond two parallel sides of the polygon passing through
the point $a_i$ and lying on the support lines to the section
$(A\cap N)$.}

\remark{\tcr Remark}{\tc To a three-dimensional set $A\subset \Bbb
RP^3$, $L$-convex-concave with respect to a line $L$, both
surgeries are applicable, since $\dim L=1=n-2$ for $n=3$. Let
$[a,b]$ be a segment on the line $L$, and  $[c,d]$ be a segment on
the line $\Bbb RP^3/L$. Then, as can be easily proved, the
surgeries $P_{[a,b]}$ and $S_{[c,d]}$ commute.}
\endremark

\subhead{\tcb A space intersecting support half-planes to
sections}
\endsubhead
{\tc Let, as before,  $A$ be a $L$-convex-concave subset of $\Bbb
RP^n$, and $\dim L=1$. Consider the following problem. Suppose
that a certain set $\{N_\alpha\}$, $\alpha \in I$, of
two-dimensional planes containing the line $L$, is fixed, and
suppose that on each affine plane $N_\alpha \setminus L$ some
supporting to a convex section $N_\alpha \cap A$ half-plane
$p^{+}_\alpha \subset N_\alpha$ is fixed. We want to find an
$(n-2)$-dimensional subspace of $\Bbb RP^n$, intersecting all
half-planes $p^{+}_\alpha$, $\alpha \in I$.}

\proclaim{\tcr Theorem}{\tci Suppose that the set $Q=\{\partial
p^+_\alpha \cap L\}$, $\alpha \in I$, contains at most four
points, where $\partial p^+_\alpha $ denotes the boundary line of
half-plane $p^{+}_\alpha$ supporting to the section  $N_\alpha
\cap A$ of a $L$-convex-concave set $A\subset \Bbb RP^n$. Then
there exists an $(n-2)$-dimensional subspace of $\Bbb RP^n$
intersecting all supporting half-planes $p^{+}_\alpha$, $\alpha
\in I$.}
\endproclaim

\demo{\tcr Proof}{\tc Suppose that the set $Q$ contains exactly
four points $a_1,\dots,a_4$ (if not, add to  $Q$ a  necessary
amount of some other points). The points $a_i$ divide the
projective line into four segments $\langle a_1, a_2\rangle$,
$\langle a_2, a_3\rangle$, $\langle a_3, a_4\rangle$, $\langle
a_4, a_1\rangle$. Denote by $I_1,\dots,I_4$ the complementary
segments  (these segments intersect one another). Apply to the set
$A$ the four surgeries $P_{I_i}$ and denote the resulting set by
$D$.

By the very definition of the set $D$ the half-planes
$p^{+}_\alpha \subset N_\alpha $ are supporting half-planes for
the sections $D\cap N_\alpha$, so any space  lying inside $D$,
will intersect the half-planes $p^{+}_\alpha$. According to the
theorem of \S5 there exists an $(n-2)$-dimensional subspace of
$\Bbb RP^n$, lying inside the set $D$. This space intersects all
half-planes $\pi^{+}_\alpha$.}
\enddemo

\enddocument